\documentclass{amsart}
\usepackage[dvips]{epsfig}
\usepackage{graphicx}
\usepackage{latexsym}
\usepackage{amsmath}
\usepackage{amsthm}
\usepackage{amssymb}
\usepackage{mathrsfs}
\usepackage[shortlabels]{enumitem}
\newlist{propenum}{enumerate}{1} 
\setlist[propenum]{label=(\roman*)}
\RequirePackage[colorlinks=true,citecolor=blue,urlcolor=blue]{hyperref}

\usepackage{fge}
\newcommand{\mysetminus}{\mathbin{\fgebackslash}}

\usepackage{caption}
\usepackage{subcaption}

\usepackage{eucal}

\usepackage{xcolor,soul}

\usepackage{mathtools}

 \usepackage[applemac]{inputenc}

\newtheorem{thm}{Theorem}[section]
\newtheorem{theo}[thm]{Theorem}

\newtheorem{lem}[thm]{Lemma}

\newtheorem{defi}[thm]{Definition}
\newtheorem{hyp}[thm]{Assumption}

\theoremstyle{remark}

\newtheorem{rem}[thm]{Remark}

\newcommand{\inL} {\underset{n \rightarrow \infty}{\overset{\rm (d)}{\xrightarrow{\hspace*{0.75cm}}}} }

\newcommand{\vip}{\vskip.2cm}
\newcommand{\COMMENTAIRE}[1]{}

\newcommand{\field}[1]{\mathbb{#1}}
\newcommand{\EE}{\field{E}}

\newcommand{\GG}{\field{G}}

\newcommand{\NN}{\field{N}}

\newcommand{\RR}{\field{R}}
\newcommand{\TT}{\field{T}}

\newcommand{\XX}{\field{X}}

\newcommand{\Bb}{{\mathcal B}}

\newcommand{\Ff}{{\mathcal F}}

\newcommand{\Hh}{{\mathcal H}}
\newcommand{\Ll}{{\mathcal L}}

\newcommand{\Nn}{{\mathcal N}}
\newcommand{\Pp}{{\mathcal P}}

\newcommand{\Qq}{{\mathcal Q}}
\newcommand{\vt}{{\vartriangle}}

\def \ep {\varepsilon}

\newcommand{\rd}{{\rm d}}

\newcommand{\cb}{{\mathcal B}}
\newcommand{\cc}{{\mathcal C}}

\newcommand{\cf}{{\mathcal F}}

\newcommand{\cp}{{\mathcal P}}
\newcommand{\cq}{{\mathcal Q}}

\newcommand{\cs}{{\mathscr S}}

\newcommand{\A}{{\mathbb A}}

\newcommand{\E}{{\mathbb E}}

\newcommand{\G}{\mathbb{G}}
\newcommand{\N}{{\mathbb N}}
\renewcommand{\P}{{\mathbb P}}

\newcommand{\R}{{\mathbb R}}

\newcommand{\T}{\mathbb{T}}

\newcommand{\ind}{{\bf 1}}

\newcommand{\sot}{\otimes_{\rm sym}}

\newcommand{\norm}[1]{\mathop{\parallel\! #1 \! \parallel}\nolimits}

\newcommand{\inv}[1]{\mathop{\frac{1}{ #1}}\nolimits}
\newcommand{\expp}[1]{\mathop {\mathrm{e}^{ #1}}}

\newcommand{\reff}[1]{(\ref{#1})}

\DeclareMathOperator*{\argmin}{arg\,min}

\newcommand{\bigO}{\mathcal{O}}

\newcommand\smallO{
  \mathchoice
    {{\scriptstyle\mathcal{O}}}
    {{\scriptstyle\mathcal{O}}}
    {{\scriptscriptstyle\mathcal{O}}}
    {\scalebox{.7}{$\scriptscriptstyle\mathcal{O}$}}
  }

\begin{document}

\title[]{Kernel estimation of the transition density in bifurcating Markov chains.}

\author{S. Val\`ere Bitseki Penda}

\address{S. Val\`ere Bitseki Penda, IMB, CNRS-UMR 5584, Universit\'e Bourgogne Franche-Comt\'e, 9 avenue Alain Savary, 21078 Dijon Cedex, France.}

\email{simeon-valere.bitseki-penda@u-bourgogne.fr}

\begin{abstract}
We study the kernel estimator of the transition density of bifurcating Markov chains. Under some ergodic and regularity properties, we prove that this estimator is consistent and asymptotically normal. Next, in the numerical studies, we propose two data-driven methods to choose the bandwidth parameters. These methods are based on the so-called two bandwidths approach.  
\end{abstract}

\maketitle



\section{Introduction}

This article is devoted to the study of the kernel estimators of the transition probability of bifurcating Markov chains. Before defining these estimators, let us first introduce useful definitions, notations and assumptions.  

\subsection{Bifurcating Markov chains}
Let $d \geq 1$ be a natural integer. In order to simplify the notations in the sequel, we set $S = \RR^{d}$ and we equip S with its Borel $\sigma\text{-}$algebra that we denote by $\cs$. We denote by $\cb(S)$  (resp. $\cb_b(S)$, resp.    $\cb_+(S)$) the  set   of  (resp.    bounded,  resp. non-negative)  $\R$-valued measurable  functions  defined  on $S$. For $f\in \cb(S)$, we set $\norm{f}_\infty =\sup\{|f(x)|, \, x\in S\}$.  For a finite measure $\lambda$ on $(S,\cs)$ and $f\in \cb(S)$ we shall write $\langle \lambda,f  \rangle$ for  $\int f(x) \,  \rd\lambda(x)$ whenever this  integral is  well defined. We denote by $\cc_b(S)$ (resp. $\cc_+(S)$) the set of bounded (resp.  non-negative) $\R$-valued continuous functions defined on $S$. For all natural integer $q \geq 1,$ we equip $S^{q}$ with $\cs^{\otimes q} = \cs \otimes \ldots \otimes \cs$, the usual product $\sigma\text{-}$field on $S^{q}$.   

Let $Q$ be a probability kernel   on $S \times \cs$, that is: $Q(\cdot  , A)$  is measurable  for all  $A\in \cs$,  and $Q(x,\cdot)$ is  a probability measure on $(S,\cs)$ for all $x \in S$. For any $f\in \cb_b(S)$,   we set for $x\in S$:
\begin{equation}\label{eq:Qf}
(Qf)(x)=\int_{S} f(y)\; Q(x,\rd y).
\end{equation}
We define $(Qf)$, or simply $Qf$, for $f\in \cb(S)$ as soon as the integral \reff{eq:Qf} is well defined, and we have $\cq f\in \cb(S)$. For $n\in \N$, we denote by $Q^n$  the $n$-th iterate of $Q$ defined by $Q^0=I_d$, the identity map on $\cb(S)$, and $Q^{n+1}f=Q^n(Qf)$ for $f\in \cb_b(S)$.  

Let $P$ be a probability kernel   on $S \times \cs^{\otimes 2}$, that is: $P(\cdot  , A)$  is measurable  for all  $A\in \cs^{\otimes 2}$,  and $P(x,\cdot)$ is  a probability measure on $(S^2,\cs^{\otimes 2})$ for all $x \in S$. For any $g\in \cb_b(S^3)$ and $h\in \cb_b(S^2)$,   we set for $x\in S$:
\begin{equation}\label{eq:Pg}
(Pg)(x)=\int_{S^2} g(x,y,z)\; P(x,\rd y,\rd z) \quad\text{and}\quad (Ph)(x)=\int_{S^2} h(y,z)\; P(x,\rd y,\rd z).
\end{equation}
We define $(Pg)$ (resp. $(Ph)$), or simply $Pg$ for $g\in \cb(S^3)$(resp. $Ph$ for $h\in \cb(S^2)$), as soon as the corresponding integral \reff{eq:Pg} is well defined, and we have  that $Pg$ and $Ph$ belong to $\cb(S)$.
\medskip 

We  now  introduce   some  notations  related  to   the  regular  binary tree.   Recall    that  $\N$ is the set of non-negative integers and $\N^*=    \N   \setminus   \{0\}$. We set  $\T_0=\G_0=\{\emptyset\}$, $\G_k=\{0,1\}^k$ and $\T_k = \bigcup  _{0   \leq  r  \leq  k}  \G_r$  for   $k\in  \N^*$,  and $\T  = \bigcup  _{r\in \N}  \G_r$.  The  set $\G_k$  corresponds to  the $k$-th generation,  $\T_k$ to the tree  up the $k$-th  generation, and $\T$ the complete binary tree. For $i\in \T$, we denote by $|i|$  the generation of $i$ ($|i|=k$ if and only if $i\in \G_k$) and $iA=\{ij; j\in A\}$ for $A\subset \T$, where $ij$ is the concatenation of the two sequences $i,j\in \T$, with the convention that $\emptyset i=i\emptyset=i$.

We recall the definition of bifurcating Markov chain (BMC) from Guyon \cite{Guyon}. 
\begin{defi}
We say  a stochastic process indexed  by $\T$, $X=(X_i,  i\in \T)$, is a bifurcating Markov chain on a measurable space $(S, \cs)$ with initial probability distribution  $\nu$ on $(S, \cs)$ and probability kernel $\cp$ on $S\times \cs^{\otimes 2}$, a  BMC in short, if:
\begin{itemize}
\item[-] (Initial  distribution.) The  random variable  $X_\emptyset$ is distributed as $\nu$.
\item[-] (Branching Markov property.) For  a sequence   $(g_i, i\in \T)$ of functions belonging to $\cb_b(S^3)$, we have for all $k\geq 0$, \[ \E\Big[\prod_{i\in \G_k} g_i(X_i,X_{i0},X_{i1}) |\sigma(X_j; j\in \T_k)\Big] =\prod_{i\in \G_k} \cp g_i(X_{i}).\]
\end{itemize}
\end{defi}

We define  three probability kernels $P_0, P_1$ and $\cq$ on $S\times \cs$ by:
\[ P_0(x,A)=\cp(x, A\times S), \quad P_1(x,A)=\cp(x, S\times A) \quad \text{for $(x,A)\in S\times  \cs$, and} \quad \cq=\inv{2}(P_0+P_1).\] 
Notice  that  $P_0$ (resp.   $P_1$)  is  the  restriction of  the  first (resp. second) marginal of $\cp$ to $S$.  Following Guyon \cite{Guyon}, we introduce an  auxiliary Markov  chain $Y=(Y_n, n\in  \N) $  on $(S,\cs)$ with  $Y_0$ distributed  as $X_\emptyset$  and transition  kernel $\cq$. The  distribution of  $Y_n$ corresponds  to the  distribution of  $X_I$, where $I$  is chosen independently from  $X$ and uniformly at  random in generation  $\G_n$.    We  shall   write  $\E_x$   when  $X_\emptyset=x$ (\textit{i.e.}  the initial  distribution  $\nu$ is  the  Dirac mass  at $x\in S$).

\medskip

Let $i,j\in \T$. We write $i\preccurlyeq  j$ if $j\in i\T$. We denote by $i\wedge j$  the most recent  common ancestor of  $i$ and $j$,  which is defined  as   the  only   $u\in  \T$   such  that   if  $v\in   \T$  and $ v\preccurlyeq i$, $v \preccurlyeq j$  then $v \preccurlyeq u$. We also define the lexicographic order $i\leq j$ if either $i \preccurlyeq j$ or $v0  \preccurlyeq i$  and $v1  \preccurlyeq j$  for $v=i\wedge  j$.  Let $X=(X_i, i\in  \T)$ be  a $BMC$  with kernel  $\cp$ and  initial measure $\nu$. For $i\in \T$, we define the $\sigma$-field:
\begin{equation*}\label{eq:field-Fi}
\cf_{i} = \sigma(X_u; u\in \T \text{ such that  $u\leq i$}).
\end{equation*}
By construction,  the $\sigma$-fields $(\cf_{i}; \, i\in \T)$ are nested as $\cf_{i}\subset \cf_{j} $ for $i\leq  j$.

For $i \in \TT$ and $k \in \NN$, we also define the $\sigma\text{-}$field:
\begin{equation*}
\Hh_{i,k} = \sigma(X_{u}, u \in i\TT_{k})
\end{equation*}

\medskip

We end this section with a useful notations.  By convention, for  $f,g\in  \cb(S)$,   we  define  the  function   $f\otimes  g$  by $(f\otimes g)(x,y)=f(x)g(y)$ for  $x,y\in S$ and 
\[
f\sot g= \inv{2}(f\otimes g + g\otimes f) \quad\text{and}\quad f\otimes ^2= f\otimes f.
\]
Notice that $\cp(g\sot \ind)=\cq(g)$ for  $g\in \cb_+(S)$.   
 
For all $u\in \TT$, we denote by $X_{u}^{\vartriangle} = (X_{u},X_{u0},X_{u1})$ the mother-daughters triangle. For a finite subset $A \subset \TT$, we define:
\begin{equation*}
M_{A}(f) = \sum_{u\in A}f(X_{u}) \quad \text{if $f \in \Bb(S)$} \quad \text{and} \quad M_{A}(f) = \sum_{u\in A}f(X_{u}^{\vartriangle}) \quad \text{if $f \in \Bb(S^{3})$}.
\end{equation*}
In the sequel we will also use the following notation: let $g$ and $h$ be two functions which depend on one variable, $x$ say; we denote by $g\oplus h$ the function of three variables, $xx_{0}x_{1} := (x,x_{0},x_{1})$ say, defined by
\begin{equation*} 
(g\oplus h)(xx_{0}x_{1}) = g(x_{0}) + h(x_{1}).
\end{equation*}

\subsection{Assumptions on the law of the bifurcating Markov chains $(X_{i}, i\in\TT)$}

For  a  set  $F\subset  \cb(S)$   of  $\R$-valued  functions,  we  write
$F^2=\{f^2; f\in F\} $, $F\otimes  F=\{f_0\otimes f_1; f_0, f_1\in F\}$,
and  $P(F)=\{Pf;  f\in F\}$  whenever  a  kernel $P$  act  on $F$.    
Following \cite{Guyon}, we state a structural assumption
on the set of functions we shall consider.

\begin{hyp}
   \label{hyp:F}
Let $F\subset \cb(S)$ be a set of $\R$-valued functions such that:
\begin{itemize}
   \item[$(i)$] $F$ is a  vector subspace which contains the constants;
   \item[$(ii)$] $F^2 \subset F$;
   \item[$(iii)$] $F\subset L^1(\nu)$; 
 \item[$(iv)$]  $F\otimes F \subset L^1(\cp(x, \cdot))$ for all $x\in S$,
    and $\cp(F\otimes F)\subset F$.
\end{itemize}
\end{hyp}
The   condition   $(iv)$   implies   that   $P_0(F)\subset   F$,$P_1(F)\subset F$ as  well as $\cq(F)\subset F$.  Notice that if  $f\in F$, then even if $|f|$ does  not belong to $F$, using  conditions  $(i)$  and $(ii)$, we get, with $g=(1+f^2)/2$, that   $|f|\leq  g$ and $g\in F$. Typically, the set $F$  can be the set  $\cc_b(S)$ of bounded real-valued functions, or the set of smooth real-valued functions such that all derivatives have  at most polynomials growth. 
\medskip

Following \cite{Guyon}, we also  consider the following  ergodic properties for $\cq$. 
\begin{hyp}\label{hyp:F1}
There exists a probability measure $\mu$ on $(S, \cs)$ such that $F\subset L^1(\mu)$ and for all $f\in F$, we have the point-wise convergence  $\lim_{n\rightarrow \infty } \cq^{n}f = \langle \mu, f \rangle$ and 
there exists $g\in F$ with:
\begin{equation}\label{eq:erg-bd}
|\cq^n(f)|\leq  g\quad\text{for all $n\in \N$.}
\end{equation}
Moreover, there exists a function $V:[1,+\infty) \mapsto (0,\infty)$ such that $V \in F$ and constants $\alpha \in (0,1)$ and $M < \infty$ such that: 
\begin{equation}\label{eq:geom-erg}
\sup_{|f| \leq V}|\cq^{n}f - \langle \mu, f \rangle| \leq M \alpha^{n} V \quad \text{for all  $n\in \N$.}
\end{equation}  
\end{hyp}
  
\begin{rem}
In particular, \eqref{eq:geom-erg} implies that for all  $f \in \Bb_{b}(S)$, we have
\begin{equation}\label{eq:geom-ergB}
|\cq^{n}f - \langle \mu, f \rangle| \leq M \|f\|_{\infty} \, \alpha^{n} \, V \quad \text{for all  $n\in \N$.}
\end{equation}
\end{rem}
 
Next, we have the following assumption on the existence of the density of $\cp.$ 
\begin{hyp}\label{hyp:DenMu}
The transition kernel $\Pp$ has a density, still denoted by $\Pp$, with respect to the Lebesgue measure. 
\end{hyp}
\begin{rem}\label{rem:}
Assumption \ref{hyp:DenMu} implies that the transition kernel $\Qq$ has a density, still denoted by $\Qq$, with respect to  the Lebesgue measure. More precisely, we have $\Qq(x,y) = 2^{-1} \int_{S} (\Pp(x,y,z)+\Pp(x,z,y))dz.$ This implies in particular that the invariant probability $\mu$ has a density, still denoted by $\mu$, with respect to the Lebesgue measure (for more details, we refer for e.g. to \cite{duflo2013random}, chap 6).
\end{rem}
 
\begin{rem}\label{rem:mu-delta}
Under Assumption \ref{hyp:DenMu}, the probability measure $\mu^{\vt}$ defined on $S^{3}$ by 
\begin{equation*}\label{eq:mu-delta}
\mu^{\vt}(dxx_{0}x_{1}) =  \mu(dx) P(x,dx_{0},dx_{1}),
\end{equation*}
has density with respect to the Lebesgue measure, that we also denote by $\mu^{\vt}$, given by $\mu^{\vt}(xx_{0}x_{1}) = \mu(x) \Pp(x,x_{0},x_{1})$, for all $xx_{0}x_{1} \in S^{3}$.
\end{rem}

\begin{hyp}\label{hyp:ub-density}
We assume that the following constant is finite:
\begin{equation*}
C_{0} = \sup_{x,x_{0},x_{1} \in S} (\mu(x) + \Qq(x,x_{0}) + \Pp(x,x_{0},x_{1})).
\end{equation*}
\end{hyp}

\begin{rem}
We recall  from \cite[Theorem~11 and Corollary~15]{Guyon}  that under Assumptions \ref{hyp:F} and \ref{hyp:F1}, we have for $f\in F$ the following convergence in $L^2$ (resp. a.s.):
\begin{equation}\label{eq:lfgn-G}
 \lim_{n\rightarrow\infty } |\G_n|^{-1} M_{\G_n}(f)=\langle \mu, f \rangle \quad \text{and} \quad \lim_{n\rightarrow\infty } |\T_n|^{-1} M_{\T_n}(f)=\langle \mu, f\rangle.
\end{equation}
\end{rem}

Now, the rest of the paper is organized follows. In Section \ref{sec:estimators}, we define the estimators of the transition density $\Pp$ based on the observation of a subpopulation. We will see that these are quotient estimators. In Section \ref{sec:asym-muvt}, we study the consistency and the asymptotic normality of the numerators of the estimators of $\Pp.$ Section \ref{sec:asym-Phat-T} is dedicated to the study of consistency and asymptotic normality of the estimators of $\Pp.$ In Section \ref{sec:numerical}, we will illustrate the consistency of our estimators in a bifurcating Markov model called bifurcating autoregressive process (BAR, for short). In particular, we will develop two data-driven bandwidth selection methods: the least squares Cross-Validation in Section \ref{sec:numerical1} and the rule of thumb type method in Section \ref{sec:numerical2}. Sections \ref{proof:thm-flx-T}-\ref{sec:cltphat-sub} are dedicated to the proofs of the main Theorems. In Section \ref{sec:proof-amise}, we prove a useful inequality and in Section \ref{sec:appendix}, we recall some useful results.  
 
\section{Kernel estimators of the transition density $\Pp$}\label{sec:estimators}
Recall that $S = \RR^{d}.$ Our aim is to estimate the transition density $\Pp$ from the observation of the subpopulation $\A_{n} \in \{\GG_{n}, \TT_{n}$\}. For that purpose, assume we observe $\XX^{\vt n} = (X^{\vt}_u)_{u \in \A_{n}}$ {\it i.e.} we have $2^{n+2} -1$ (or $3\times2^{n}$) random variables with value in $S$. Let $K_{0}: S \rightarrow \RR$ and $K: S^{3} \rightarrow \RR$ be a functions such that $\int_{S} K_{0}(x) dx = 1$ and $K = K_{0}\otimes K_{0} \otimes K_{0}$. We also have $\int_{S^{3}} K(xx_{0}x_{1})dxx_{0}x_{1} = 1.$ Let $(h_{n},n \in \NN)$ be a sequence of positive numbers which converges to $0$ as $n$ goes to infinity. When there is no ambiguity, we write $h$ for $h_{n}$. Let $\A_{n} \in \{\TT_{n}, \GG_{n}\}$. We define, for all $x \in S$: 
\begin{equation}\label{eq:def-estim-mu}
\widehat{\mu}_{\A_{n}}(x) = \frac{1}{|\A_{n}|h^{d/2}} \sum_{u \in \A_{n}} K_{0h_{n}}(x - X_{u}), 
\end{equation}
where $K_{0h_{n}}(x-y) = h_{n}^{-d/2}K_{0}(h_{n}^{-1}(x-y))$ and for all $xx_{0}x_{1} \in S^{3}$:
\begin{equation}\label{eq:def-estim}
\widehat{\mu}^{\vartriangle}_{\A_{n}}(xx_{0}x_{1}) = \frac{1}{|\A_{n}|h^{3d/2}}\sum_{u \in \A_{n}} K_{h_{n}}(xx_{0}x_{1} - X_{u}^{\vartriangle}) \quad \text{and} \quad \widehat{\Pp}_{\A_{n}}(xx_{0}x_{1}) = \frac{\widehat{\mu}^{\vartriangle}_{\A_{n}}(xx_{0}x_{1})}{\widehat{\mu}_{\A_{n}}(x)},
\end{equation}
where
\begin{equation*} 
K_{h_{n}}(xx_{0}x_{1} - yy_{0}y_{1}) =  h_{n}^{-3d/2} \, K\left(h_{n}^{-1}(x-y), h_{n}^{-1}(x_{0} - y_{0}), h_{n}^{-1}(x_{1} - y_{1})\right), 
\end{equation*} 
with the convention that $\widehat{\Pp}_{\A_{n}}(xx_{0}x_{1}) = 0$ if $\widehat{\mu}_{\A_{n}}(x) = 0.$ However, we stress that if we assume that $K_{0}$ is strictly positive, then $\widehat{\mu}_{\A_{n}}(x) > 0$ for all $x \in S.$

\medskip

From now on, we fix $xx_{0}x_{1} \in S^{3}$, that is, we are interested in the estimation at the point $xx_{0}x_{1}$. We assume that $\mu(x) \neq 0.$ We consider the function $f_{n}$ defined by:
\begin{equation}\label{eq:def-fnT}
f_{n}(yy_{0}y_{1}) = K_{h}(xx_{0}x_{1} - yy_{0}y_{1}).  
\end{equation}
If we want to be more rigorous, we must write $f_{n,xx_{0}x_{1}}$ instead of $f_{n}$. But, we choose to write without the index $xx_{0}x_{1}$ in order to simplify the writing.

\begin{rem}
Note that asymptotic behavior (consistence and asymptotic normality) of $\widehat{\mu}_{\A_{n}}$ have been studied in \cite{BD2020}.
\end{rem}

\begin{rem}\label{rem:ext-h}
We stress that the results of this paper can be straightforward extended the case where the bandwidth $\pmb{h}$ is a vector of $\RR^{3d}$, with possibly different coordinates. More precisely, one can take the bandwidth $\pmb{h} = (h_{i}, 1\leq i \leq 3d),$ where the $h_{i}$'s may take different values. For our convenience, we choose to work with the case where all the coordinates are the same, that is $h_{i} = h$ for all $1 \leq i \leq 3d.$
\end{rem}

\section{Consistency and Asymptotic normality for $\widehat{\mu}^{\vt}_{\A_{n}}(xx_{0}x_{1})$}\label{sec:asym-muvt}

First, we will study the consistency and the asymptotic normality of $\widehat{\mu}^{\vt}_{\A_{n}}(xx_{0}x_{1}).$ We set $\tilde{f}_{n} = f_{n} - \langle \mu,\Pp f_{n} \rangle$. We begin with the study asymptotic normality of $N_{n,\emptyset}(f_{n}) = |\GG_{n}|^{-1/2} M_{\A_{n}}(\tilde{f}_{n}).$ This is motivated by the following decomposition:
\begin{equation}\label{eq:Dmuhatvt-T}
\widehat{\mu}^{\vartriangle}_{\A_{n}}(xx_{0}x_{1}) - \mu^{\vt}(xx_{0}x_{1}) = (|\A_{n}| \, |\GG_{n}|^{-1/2} \, h^{3d/2})^{-1} N_{n,\emptyset}(f_{n}) \, + \, (h^{-3d/2} \langle \mu^{\vt}, f_{n}\rangle - \mu^{\vt}(xx_{0}x_{1})).
\end{equation}
We will need the following assumption on the bandwidth and on the kernel. 
\begin{hyp}\label{hyp:2alpha2h3d}
$ $

We assume that:
\begin{itemize}
\item[(i)] $h_{n} = 2^{-n \gamma}$ and $2\alpha^{2} < 2^{3d\gamma}$ for some $\gamma \in (0,1/3d)$.  
\item[(ii)] The kernel $K_{0}$ (resp. $K_{0}^{2}$) is integrable and square integrable.
\end{itemize}
\end{hyp}
\begin{rem}\label{rem:2alpha2h3d}
Assumption \ref{hyp:2alpha2h3d}, (i) implies in particular that 
\begin{equation}\label{eq:2alpha2h3d}
\lim_{n \rightarrow \infty} |\GG_{n}|h_{n}^{3d} = \infty \quad  \text{and} \quad \lim_{n \rightarrow \infty}(2\alpha^{2})^{n} \, h_{n}^{3d} = 0.
\end{equation}
Note that Assumption \ref{hyp:2alpha2h3d}, (i) is automatically satisfied if $2\alpha^{2} \leq 1$, regardless of the value of $\alpha$. For $2\alpha^{2} > 1$, this Assumption implies that the choice of the bandwidth is function of the ergodicity rate of the auxiliary Markov chain $Y$.
\end{rem}

We have the following result.

\begin{thm}\label{thm:flx-T}
Let  $X$  be  a  BMC   with  kernel  $\cp$  and  initial  distribution $\nu$ such that Assumptions \ref{hyp:F},  \ref{hyp:F1}, \ref{hyp:DenMu}, \ref{hyp:ub-density} and \ref{hyp:2alpha2h3d} hold. Then, we  have the  following convergence  in distribution: 
\begin{equation*}
N_{n, \emptyset}(f_{n}) \; \xrightarrow[n\rightarrow \infty ]{\text{(d)}} \; G,
\end{equation*}
where $G$ is a centered Gaussian  random variable with finite variance $\sigma^2 = 2\, \|K_{0}\|_{2}^{6}\, \mu^{\vt}(x,x_{0},x_{1})$ if $\A_{n} = \TT_{n}$ and $\sigma^2 = \|K_{0}\|_{2}^{6}\, \mu^{\vt}(x,x_{0},x_{1})$ if $\A_{n} = \GG_{n}$. 
\end{thm}
\begin{proof}
The proof of Theorem \ref{thm:flx-T} is postponed to Section \ref{proof:thm-flx-T}.
\end{proof}
Next, in order to study the asymptotic normality of $\widehat{\mu}^{\vt}_{\A_{n}}(xx_{0}x_{1}),$ we do the following additional hypothesis.

\begin{hyp}\label{hyp:estim-tcl}
We assume that Assumption \ref{hyp:2alpha2h3d} holds and there exists  $s > 0$ such that the following holds. 
\begin{itemize}
\item[(iv)]\textbf{The density $\mu^{\vt}$ (resp. $\mu$) belongs to the (isotropic) H\"older class of order $(s, \ldots, s) \in \RR^{3d}$ (resp. $(s, \ldots, s) \in \RR^{d}$):}  The density  $\mu^{\vt}$ admits partial derivatives with respect to $x_{j}$, for all $j\in \{1, \ldots 3d\}$,  up to the order $\lfloor s \rfloor$ and there exists a finite constant  $L > 0$ such that  for all $x=(x_1, \ldots, x_{3d}), \in \RR^{3d}$, $t\in \R$ and $ j \in \{1, \ldots, 3d\}$: 
\begin{equation*}
\left|\frac{\partial^{\lfloor s \rfloor}\mu^{\vt}}{\partial x_{j}^{\lfloor s \rfloor}}(x_{-j},t)-\frac{\partial^{\lfloor s \rfloor}\mu^{\vt}}{\partial x_{j}^{\lfloor s \rfloor}}(x)\right| \leq L|x_{j} - t|^{\{s\}}, 
\end{equation*}
where $(x_{-j},t)$ denotes the vector $x$ where we have replaced the $j^{th}$ coordinate $x_{j}$ by $t$, with the convention ${\partial^{0}\mu^{\vt}}/{\partial x_{j}^{0}} = \mu^{\vt}$. The same thing for the density $\mu.$
\item[(v)]\textbf{The kernel $K_{0}$ is  of order $(\lfloor s \rfloor, \ldots, \lfloor s \rfloor) \in \NN^{d}$:} We have $\int_{\RR^{d}} |x|^{s}K_{0}(x)\, dx < \infty$  and $\int_{\RR} x^{k}_{j}\, K_{0}(x)\, dx_{j} = 0$ for all $k \in \{1,\ldots,\lfloor s \rfloor\}$ and $j \in \{1,\ldots,d\}$.
\item[(vi)]\textbf{Bandwith control:} We have  $\gamma> 1/(2s+3d)$, that is $\lim_{n \rightarrow \infty} |\GG_{n}|h_{n}^{2s + 3d} = 0$. 
\end{itemize}
\end{hyp}
Notice that Assumption \ref{hyp:estim-tcl}-$(iv)$ implies that $\mu^{\vt}$ (resp. $\mu$) is at least  H\"older continuous as $s>0$. We have the following result.

\begin{thm}\label{thm:cltEst-T}
Let  $X$  be  a  BMC   with  kernel  $\cp$  and  initial  distribution $\nu$. Under Assumptions of Theorem \ref{thm:flx-T}, we have for all $(x,x_{0},x_{1})$ 
and $\A_{n} \in \{\GG_{n}, \TT_{n}\}$
\begin{equation}
\widehat{\mu}_{\A_{n}}^{\vt}(xx_{0}x_{1}) \xrightarrow[n\rightarrow \infty ]{\P}  \mu^{\vt}(xx_{0}x_{1}) \quad\text{in probability}. \label{eq:cge-mu-tri}
\end{equation}
Moreover, under the additional Assumption \ref{hyp:estim-tcl}, we  have 
the  following convergence in distribution: 
\begin{equation*}
|\A_{n}|^{1/2} \, h_{n}^{3d/2} \, (\widehat{\mu}_{\A_{n}}^{\vt}(xx_{0}x_{1}) - \mu^{\vt}(xx_{0}x_{1})) \; \xrightarrow[n\rightarrow \infty ]{\text{(d)}} \; G \quad\text{in distribution}, \nonumber
\end{equation*}
where $G$ is a centered Gaussian  random variable with finite variance $\sigma^2 = \|K_{0}\|_{2}^{6}\, \mu^{\vt}(x,x_{0},x_{1})$. 
\end{thm}
\begin{proof}
The proof is postponed to Section \ref{proof:thm:cltEst-T}.
\end{proof}

\section{Consistency and Asymptotic normality for $\widehat{\Pp}^{\vt}_{\A_{n}}(xx_{0}x_{1})$}\label{sec:asym-Phat-T}

We are now in position to state consistency and asymptotic normality of kernel estimator of the transition density $\Pp$.
First, as a consequence of \eqref{eq:def-estim}, \eqref{eq:cge-mu-tri} and \eqref{eq:LSTnmuhat-mu} below, we have the following result.

\begin{lem}
Under the Assumptions of Theorem \ref{thm:flx-T}, we have for all $(x,x_{0},x_{1})$ 
and $\A_{n} \in \{\GG_{n}, \TT_{n}\}:$
\[
\widehat{\Pp}_{\A_{n}}(x,x_{0},x_{1}) \xrightarrow[n\rightarrow \infty ]{\P}  \Pp(x,x_{0},x_{1}) \quad\text{in probability}.
\]
\end{lem} 
Next, we have the following result.
\begin{theo}\label{thm:cltPhat-sub}
Let  $X$  be  a  BMC   with  kernel  $\cp$  and  initial  distribution $\nu$. Under the assumptions of Theorem \ref{thm:flx-T} and the additional Assumption \ref{hyp:estim-tcl}, we have, 
\begin{equation*}\label{eq:Phat-sub}
\sqrt{|\A_{n}|h_{n}^{3d}} (\widehat{\Pp}_{\A_{n}}(x,x_{0},x_{1}) - \Pp(x,x_{0},x_{1})) \inL G,
\end{equation*}
where $G$ is a centered Gaussian real-valued random variable with mean $0$ and variance $$\sigma^{2} \, = \, \|K_{0}\|_{2}^{6}\, \Pp(x,x_{0},x_{1})/\mu(x).$$
\end{theo}

\begin{proof}
The proof is postponed to Section \ref{sec:cltphat-sub}.
\end{proof}

\section{Numerical studies}\label{sec:numerical}

We consider the real-valued Gaussian bifurcating autoregressive
process (BAR) $X=(X_{u},u\in\TT)$ where $X_{\emptyset}$ is arbitrary and for all $ u \in \T$: 
\begin{equation}\label{eq:bar-asym}
\begin{cases}
X_{u0} = a_{0} X_{u} + b_{0} + \ep_{u0}  \\
 X_{u1} = a_{1} X_{u} + b_{1} + \ep_{u1},
\end{cases}
\end{equation}
with $a_{0}, a_{1} \in [-1,1]$, $b_{0}, b_{1} \in \RR$ and
$((\ep_{u0},\ep_{u1}),\, u \in \T)$  an
independent  sequence of bivariate Gaussian  $\Nn(0,\Gamma)$ random vectors
independent of $X_{\emptyset}$ with covariance matrix, with $\sigma>0$
and $\rho\in \R$ such that $|\rho| < \sigma^2$:
\begin{equation*}
\Gamma = \begin{pmatrix} \sigma^{2} \quad \rho \\ \rho \quad \sigma^{2} \end{pmatrix}.
\end{equation*}
Then the process $X=(X_{u},u\in\TT)$ is a BMC with transition
probability $\Pp$ given by:
\[
\Pp(x,dy,dz) = \frac{1}{2\pi \sqrt{\sigma^{4} - \rho^{2}}}
\, \exp\left( -\frac{\sigma^{2}}{2(\sigma^{4} -
  \rho^{2})}\,  g(x,y,z) \right)\, dydz,
\]
with 
\[
g(x,y,z)=(y-a_{0}x-b_{0})^{2}  - 2\rho\sigma^{-2}(y - a_{0}x -
b_{0})(z - a_{1}x - b_{1}) + (z - a_{1}x - b_{1})^{2}.
\]
The transition kernel $\Qq$ of the auxiliary Markov  chain is defined by:
\begin{equation*}
\Qq(x,dy) = \frac{1}{2\sqrt{2\pi\sigma^2}}\left(\expp{-(y - a_{0}x
    - b_{0})^{2}/2\sigma^{2}}  + \expp{-(y - a_{1}x -
    b_{1})^{2}/2\sigma^{2}}\right)\, dy.
\end{equation*}
We will estimate the transition density $\Pp$ in a compact set $D \subset \RR^{3}$. For that purpose, we use the estimator $\widehat{\Pp}_{\GG_{n}}(xx_{0}x_{1}),$ for all $xx_{0}x_{1} \in D$, given in \eqref{eq:def-estim}, with the Gaussian kernel $K_{0}$ defined by
\begin{equation}\label{eq:GaussK}
K_{0}(x) = \frac{1}{\sqrt{2\pi}} \expp{-x^2/2}.
\end{equation}
Since the bandwidth is a function of the ergodicity rate which is unknown, we have to develop a method based on data in order to select it. To select the optimal bandwidth for $\widehat{\Pp}_{\GG_{n}}$ defined in \eqref{eq:def-estim}, we will use the so-called ``two bandwidths approach'' (see for e.g. \cite{comte-marie2021}). More precisely, since $\widehat{\Pp}_{\GG_{n}}$ is a quotient estimator, we select separately the bandwidths for the numerator ($h_{N}$, say) and the denominator ($h_{D}$, say). For that purpose, we propose two methods: the cross validation and the rule of thumb type method. The objective here is not to study nor to compare theoretically these two methods. This will be done in the future works. Our objective is only the see the numerical performances of each method. 
%
Our conclusion is that even if the rule of thumb developed in this paper give a crude approximation, it as more computational benefit with respect to the least squared cross validation.

\subsection{Bandwidth selection by least squares Cross-Validation method}\label{sec:numerical1}
We choose the bandwidths which minimises the mean integrated squared errors (MISEs)
\begin{equation*}
\EE\left[\int_{\RR^{3}} (\widehat{\mu}^{\vt}_{\GG_{n}} - \mu^{\vt})^{2}(xx_{0}x_{1})dxdx_{0}dx_{1}\right] \quad \text{and} \quad \EE\left[\int_{\RR} (\widehat{\mu}_{\GG_{n}} - \mu)^{2}(x)dx\right], 
\end{equation*}
where $\widehat{\mu}_{\GG_{n}}^{\vt}$ and $\widehat{\mu}_{\GG_{n}}$ are defined in \eqref{eq:def-estim-mu} and \eqref{eq:def-estim}. This is equivalent to minimise the functions $J^{\vt}$ and $J$ defined by
\begin{equation*}
J^{\vt}(h) = \EE\left[\int_{\RR^{3}}(\widehat{\mu}_{\GG_{n}}^{\vt})^{2}(xx_{0}x_{1})dxdx_{0}dx_{1}\right] - 2 \EE\left[ \int_{\RR^{3}} (\widehat{\mu}_{\GG_{n}}^{\vt} \mu^{\vt})(x,x_{0},x_{1}) dxdx_{0}dx_{1}\right]
\end{equation*}
and
\begin{equation*}
J(h) = \EE\left[\int_{\RR}(\widehat{\mu}_{\GG_{n}})^{2}(x)dx\right] - 2 \EE\left[ \int_{\RR} (\widehat{\mu}_{\GG_{n}} \mu)(x) dx\right].
\end{equation*}
The method to select the bandwidths is the following.
\begin{enumerate}
\item We divide the sample $(X^{\vt}_{u})_{u \in \GG_{n}}$ into $K$ disjoints subsamples $\{(X^{\vt}_{u})_{u \in \GG_{n}^{(k)}}, k \in \{1, \ldots, K\}\},$ with $(\GG_{n}^{(k)}, k \in \{1, \ldots, K\})$ a partition of $\GG_{n}.$
\item For each subsample $(X_{u})_{u \in \GG_{n}^{(k)}}$: 
\begin{enumerate}
\item We set $\widehat{\mu}^{\vt}_{\GG_{n}[-k]}$ and $\widehat{\mu}_{\GG_{n}[-k]}$ the estimators of $\mu^{\vt}$ and $\mu$ obtaining using the subsample  $(X^{\vt}_{u})_{u \in \GG_{n}} \mysetminus (X^{\vt}_{u})_{u \in \GG_{n}^{(k)}},$ where for two sets $A$ and $B$, $B \mysetminus A$ denotes the set of elements in $B$ but not in $A.$ More precisely, 
\begin{equation*}
\widehat{\mu}^{\vt}_{\GG_{n}[-k]}(xx_{0}x_{1}) = \frac{1}{|\GG_{n}[-k]|h^{3}} \sum_{u \in \GG_{n}[-k]} K_{0}\left(\frac{x-X_{u}}{h}\right)K_{0}\left(\frac{x_{0}-X_{u0}}{h}\right)K_{0}\left(\frac{x_{1}-X_{u1}}{h}\right)
\end{equation*}
and
\begin{equation*}
\widehat{\mu}_{\GG_{n}[-k]}(x) = \frac{1}{|\GG_{n}[-k]|h} \sum_{u \in \GG_{n}[-k]} K_{0}\left(\frac{x-X_{u}}{h}\right), 
\end{equation*}
where we set $\GG_{n}[-k] = \GG_{n} \mysetminus \GG_{n}^{(k)}.$
\item We approximate $J$ and $J^{\vt}$ by
\begin{align*}
&\widehat{J}^{(k)}(h) = \int_{\RR}(\widehat{\mu}_{\GG_{n}[-k]})^{2}(x)dx - \frac{2}{|\GG_{n}^{(k)}|} \sum_{u \in \GG_{n}^{(k)}} \widehat{\mu}_{\GG_{n}[-k]}(X_{u}) \quad \text{and}\\
&\widehat{J}^{\vt(k)}(h) = \int_{\RR^{3}}(\widehat{\mu}^{\vt}_{\GG_{n}[-k]})^{2}(xx_{0}x_{1})dxdx_{0}dx_{1} - \frac{2}{|\GG_{n}^{(k)}|} \sum_{u \in \GG_{n}^{(k)}} \widehat{\mu}^{\vt}_{\GG_{n}[-k]}(X_{u},X_{u0},X_{u1}).
\end{align*}
\end{enumerate}
\item Let $\Hh = \{h_{1}, \ldots, h_{m}\} \subset (0,1]$ be a bandwidth grid. Then, the selected bandwidths $\widehat{h}_{N}$ and $\widehat{h}_{D}$ for the numerator and the denominator of $\widehat{\Pp}_{\GG_{n}}$ are given by:
\begin{align*}
\widehat{h}_{N} :=  \argmin_{h \in \Hh} \frac{1}{K} \sum_{k=1}^{K} \widehat{J}^{\vt(k)}(h) \quad \text{and} \quad \widehat{h}_{D} :=  \argmin_{h \in \Hh} \frac{1}{K} \sum_{k=1}^{K} \widehat{J}^{(k)}(h).
\end{align*} 
\end{enumerate}
Finally, the estimator used for numerical studies is $\widetilde{\Pp}_{\GG_{n}}$ defined by:
\begin{equation*}
\widetilde{\Pp}_{\GG_{n}}(xx_{0}x_{1}) = \frac{\widetilde{\mu}^{\vartriangle}_{\A_{n}}(xx_{0}x_{1})}{\widetilde{\mu}_{\A_{n}}(x)},
\end{equation*}
with
\begin{align*}
&\widetilde{\mu}_{\GG_{n}}(x) = \frac{1}{|\GG_{n}|\widehat{h}_{D}} \sum_{u \in \GG_{n}} K_{0}\left(\frac{x - X_{u}}{\widehat{h}_{D}}\right)  \quad \text{and} \quad \\
&\widetilde{\mu}^{\vartriangle}_{\GG_{n}}(xx_{0}x_{1}) = \frac{1}{|\GG_{n}|\widehat{h}_{N}^{3}}\sum_{u \in \GG_{n}} K_{0}\left(\frac{x - X_{u}}{\widehat{h}_{N}}\right) K_{0}\left(\frac{x_{0} - X_{u0}}{\widehat{h}_{N}}\right) K_{0}\left(\frac{x_{1} - X_{u1}}{\widehat{h}_{N}}\right). 
\end{align*}
This method is known as the $K\text{-}$fold cross validation. One advantage of this method in the context of bifurcating Markov chains is that it is not requires the knowledge of the ergodicity rate. The main drawback being that it requires a lot of time for calculations.

\subsection{Gaussian symmetric BAR reference bandwidth selection}\label{sec:numerical2}

In order to define a selection rule, we consider the special case of Gaussian BAR defined by \eqref{eq:bar-asym} where $a_{0} = a_{1} := a$ and $\rho = 0$ as a reference model. It is well known (see \cite{BD1}) that  the densities of the transition kernel $\Qq$ of the auxiliary Markov chain and the invariant probability $\mu$ associated to $\Qq$ are given by:
\[
\Qq(x,y) = \frac{1}{\sqrt{2\pi\sigma^2}}\exp\left(-\frac{(y - a x )^{2}}{2\sigma^{2}}\right) \quad \text{and}  \quad \mu(x) = \frac{1}{\sqrt{2\pi}\sigma_{a}} \exp\left(-\frac{x^{2}}{2 \sigma_{a}^{2}}\right), 
\]
where $\sigma_{a} = \sigma/\sqrt{1 - a^{2}}.$
The density of the transition probability $\Pp$ associated to this bifurcating Markov chain is defined by $\Pp(x,y,z) = \Qq(x,y) \Qq(x,z)$ and then, $\mu^{\vt}(x,y,z) = \mu(x)\Qq(x,y) \Qq(x,z).$
We 
then have that the invariant densities $\mu$ and $\mu^{\vt}$ are square integrable and twice differentiable, with the second order derivative of $\mu$ and all the second order partial derivatives of $\mu^{\vt}$ bounded, continuous and square integrable. It is also well that the Markov chain with transition $\Qq$ is geometrically ergodic and that the geometric ergodic rate of convergence is $a$ (for more details, see for e.g Example 2.8 in \cite{BD1}). In particular, following the proof of Proposition 28 in \cite{Guyon}, one can prove that for all derivable function $f$ such $f$ and $f'$ are bounded, we have 
\begin{equation}\label{eq:erg-BAR}
|\Qq^{n}f(x) - \langle \mu, f \rangle| \leq \|f'\|_{\infty} (\sigma(1 - a)^{-1} + |x|) a^{n}.
\end{equation}
We assume that $\Ll(X_{\emptyset}) = \mu$, that is $X_{\emptyset}$ is distributed as $\mu$, which implies that the process is stationary. We are now going to behave as if we did not know the invariant measures $\mu$ and $\mu^{\vt}$ and the transition probability $\Pp$. Recall the kernel density estimator of $\mu$ defined in \eqref{eq:def-estim-mu} and the kernel $K_{0}$ defined in \eqref{eq:GaussK}. Recall also the kernel estimator of the transition density $\Pp$ defined in \eqref{eq:def-estim}. To ease notation, we write $\widehat{\mu}$ and $\widehat{\mu}^{\vt}$ instead of $\widehat{\mu}_{\GG_{n}}$ and $\widehat{\mu}_{\GG_{n}}^{\vt}$ respectively. We recall that our strategy is to select bandwidth for the numerator and the denominator in the estimation of $\Pp$. First, we treat the denominator $\widehat{\mu}.$ The selection rule is based on the following asymptotic upper bound, known as asymptotic mean squared error:
\begin{equation}\label{eq:AMSE}
\EE\left[\left( \widehat{\mu}(x) - \mu(x) \right)^{2}\right] \leq \frac{h^{4}}{2} \kappa_{2}^{2} \mu''(x)^{2} \, + \, \frac{4 \|K_{0}\|_{2}^{2}}{|\GG_{n}|h} \mu(x)  \, + \, \frac{2 \, C_{a,\sigma}}{a^{2} |\GG_{n}|} \sum_{k=1}^{n-1} (2a^{2})^{k}   \, + \,  \smallO(1).
\end{equation}
where 
\begin{equation}\label{eq:kappa2}
C_{a,\sigma} \, = \, \frac{1}{\pi \sigma^{2} (1+a)(1-a)^{2}} \quad  \text{and} \quad \kappa_{2} \, = \, \int_{\RR} y^{2} K_{0}(y) dy = 1.
\end{equation}
We postponed the proof of \eqref{eq:AMSE} in Section \ref{sec:proof-amise}. Now, let $p$ be a non negative probability density defined in $\RR$ such that $\|p\|_{\infty} \leq 1.$ Then, \eqref{eq:AMSE} implies that
\begin{equation}\label{eq:p-AMISE}
\int_{\RR} \EE\left[\left( \widehat{\mu}(x) - \mu(x) \right)^{2}\right] p(x) dx \leq \frac{h^{4}}{2} \kappa_{2}^{2} \int_{\RR} \mu''(x)^{2} dx \, + \, \frac{4 \|K_{0}\|_{2}^{2}}{|\GG_{n}|h}  \, + \, \frac{2 \, C_{a,\sigma}}{a^{2} |\GG_{n}|} \sum_{k=1}^{n-1} (2a^{2})^{k}   \, + \,  \smallO(1).
\end{equation}
The term in the left hand side of \eqref{eq:p-AMISE} is a modification of asymptotic mean integrated squared error that we call $p\text{-}$AMISE. We have introduced it because the last term in \eqref{eq:AMSE} does not depend on $x$. Finally, \eqref{eq:p-AMISE} suggests us to choose the bandwidth which minimises the function $G$ defined by
\begin{equation*}
G(h) \, = \, \frac{h^{4}}{2} \kappa_{2}^{2} \int_{\RR} \mu''(x)^{2} dx \, + \, \frac{4 \|K_{0}\|_{2}^{2}}{|\GG_{n}|h}  + M_{a,\sigma} \, a^{2n} \, \ind_{\{2a^{2} > 1\}},
\end{equation*}
where $M_{a,\sigma} = 2/(\pi \sigma^{2} a^{2} (1 + a) (1 - a)^{2} (2 a^{2} - 1)).$ Optimizing in $h$, we get that the optimal bandwidth $h_{D}$ (for the denominator of $\widehat{\Pp}_{\GG_{n}}$ defined in \eqref{eq:def-estim})  is given by
\begin{equation*}
h_{D} = (c_{1}/4c_{2})^{1/5} |\GG_{n}|^{-1/5} \, \ind_{\{2 a^{2} \leq 2^{-1/5}\}} \, + \, (c_{1}/M_{a,\sigma}) \, (2 a^{2})^{-n} \, \ind_{\{2 a^{2} > 2^{-1/5}\}},
\end{equation*}
where
\begin{equation*}
c_{1} = 4 \|K_{0}\|_{2}^{2} = \frac{2}{\sqrt{\pi}} \quad \text{and} \quad c_{2} = (1/2) \kappa_{2}^{2} \int_{\RR} \mu''(x)^{2} dx = \frac{3}{16 \sqrt{\pi}} \, \sigma_{a}^{-5}.
\end{equation*}
Next, we treat the numerator $\widehat{\mu}^{\vt}$ of $\widehat{\Pp}_{\GG_{n}}.$ Recalling Remark \ref{rem:ext-h}, we consider the general case where for all $xx_{0}x_{1} \in S^{3}$:
\begin{equation*}
\widehat{\mu}^{\vartriangle}_{\GG_{n}}(xx_{0}x_{1}) = \frac{1}{|\GG_{n}| h h_{0} h_{1}}\sum_{u \in \GG_{n}} K_{0}(h^{-1}(x-X_{u})) K_{0}(h_{0}^{-1}(x_{0}-X_{u0})) K_{0}(h_{1}^{-1}(x_{1}-X_{u1})).
\end{equation*}
Recall that for a vector $\pmb{v},$ $\pmb{v}^{t}$ denotes its transpose. As in \eqref{eq:AMSE}, we have the following asymptotic upper bound:
\begin{multline}\label{eq:AMSE-vt}
\EE\left[\left( \widehat{\mu}^{\vt}(xx_{0}x_{1}) - \mu^{\vt}(xx_{0}x_{1}) \right)^{2}\right] \, \leq \, \frac{1}{2} \kappa_{2}^{2} \, \left(\pmb{h}^{t} \, \pmb{H}(xx_{0}x_{1}) \, \pmb{h}\right)^{2} + \frac{6 \, \|K_{0}\|_{2}^{6}}{|\GG_{n}|hh_{0}h_{1}} \mu^{\vt}(xx_{0}x_{1}) \\ 
+ \, \frac{C^{\vt}_{a,\sigma} \, \Pp(xx_{0}x_{1})^{2}}{|\GG_{n}|} \, \sum_{k=1}^{n-1} (2 \, a^{2})^{k} \, + \, \smallO(1),
\end{multline}
where $\kappa_{2}$ is defined in \eqref{eq:kappa2}, $\pmb{h} = (h, h_{0}, h_{1})^{t},$ 
and
\begin{align*}
&\pmb{H}(xx_{0}x_{1}) = \begin{pmatrix} 
\frac{\partial^{2} \mu^{\vt}}{\partial x^{2}}(x,x_{0},x_{1}) & 0  & 0 \\ 0 & \frac{\partial^{2} \mu^{\vt}}{\partial x_{0}^{2}}(x,x_{0},x_{1}) & 0 \\ 0 & 0 & \frac{\partial^{2} \mu^{\vt}}{\partial x_{1}^{2}}(x,x_{0},x_{1})  
\end{pmatrix}, \\
& C^{\vt}_{a,\sigma} = \frac{4}{e \, \pi \sigma^{2} a^{2} (1 - a)^{2} (1 + a)}.
\end{align*}
We let the proof of \eqref{eq:AMSE-vt} to the reader since it follows the same lines that of \eqref{eq:AMSE}. Let $p$ be a non negative probability density defined in $\RR$ such that $\|p\|_{\infty} \leq 1.$ Integrating \eqref{eq:AMSE-vt} with respect to $p(x)dxdx_{0}dx_{1}$, we get
\begin{multline*}
\iiint_{\RR^{3}} \EE\left[\left( \widehat{\mu}^{\vt}(xx_{0}x_{1}) - \mu^{\vt}(xx_{0}x_{1}) \right)^{2}\right] p(x)dxdx_{0}dx_{1} \leq \frac{1}{2} \kappa_{2}^{2} \, \iiint_{\RR^{3}} \left(\pmb{h}^{t} \, \pmb{H}(xx_{0}x_{1}) \, \pmb{h}\right)^{2} dxdx_{0}dx_{1} \\
+ \, \frac{6 \, \|K_{0}\|_{2}^{6}}{|\GG_{n}|hh_{0}h_{1}} \, + \, \frac{C^{\vt}_{a,\sigma}}{4 \pi \sigma^{2}|\GG_{n}|} \, \sum_{k=1}^{n-1} (2 \, a^{2})^{k} \, + \, \smallO(1).  
\end{multline*}
Now, the latter equation suggests us to choose the vector bandwidth $\pmb{h}$ which minimises the function $G^{\vt}$ defined by
\begin{equation*}
G^{\vt}(h,h_{0},h_{1}) =  \frac{1}{2} \kappa_{2}^{2} \, \iiint_{\RR^{3}} \left(\pmb{h}^{t} \, \pmb{H}(xx_{0}x_{1}) \, \pmb{h}\right)^{2} dxdx_{0}dx_{1} \\
+ \, \frac{6 \, \|K_{0}\|_{2}^{6}}{|\GG_{n}|hh_{0}h_{1}} \, + \, M^{\vt}_{a,\sigma} \, a^{2n} \ind_{\{2a^{2} > 1\}}.
\end{equation*}
where $M^{\vt}_{a,\sigma} = (1/(4 \, \pi \, \sigma^{2} (2a^{2} - 1))) \, C^{\vt}_{a,\sigma}.$ Optimizing the function $G^{\vt}$ in $\pmb{h},$ we get that the optimal bandwidth $\pmb{h}_{N} = (h_{N}, h_{0N}, h_{1N})$ is given by
\begin{align*}
&h_{N} \, = \, \left(c^{\vt}_{2}/(4 c^{\vt})\right)^{1/7} |\GG_{n}|^{-1/7} \, \ind_{\{2a^{2} \leq 2^{-3/7}\}} \, + \, (c^{\vt}_{2}/M^{\vt}_{a,\sigma})^{1/3}  \, (2a^{2})^{-n/3}  \, \ind_{\{2a^{2} > 2^{-3/7}\}}, \\
&h_{0N} \, = \, \left(c^{\vt}_{2}/(4 c^{\vt}_{0})\right)^{1/7} |\GG_{n}|^{-1/7} \, \ind_{\{2a^{2} \leq 2^{-3/7}\}} \, + \, (c^{\vt}_{2}/M^{\vt}_{a,\sigma})^{1/3} \,  (2a^{2})^{-n/3}  \, \ind_{\{2a^{2} > 2^{-3/7}\}}, \\
&h_{1N} \, = \, \left(c^{\vt}_{2}/(4 c^{\vt}_{1})\right)^{1/7} |\GG_{n}|^{-1/7} \, \ind_{\{2a^{2} \leq 2^{-3/7}\}} \, + \, (c^{\vt}_{2}/M^{\vt}_{a,\sigma})^{1/3} \,   (2a^{2})^{-n/3}  \, \ind_{\{2a^{2} > 2^{-3/7}\}},
\end{align*}
where $\kappa_{2}$ is defined in \eqref{eq:kappa2}, $c^{\vt}_{2} = 6 \, \|K_{0}\|_{2}^{6} = 6/(8 \pi \sqrt{\pi})$ and
\begin{align*}
&c^{\vt} \, =  \, \frac{1}{2} \, \kappa_{2}^{2} \iiint_{\RR^{3}} \left(\frac{\partial^{2} \mu^{\vt}}{\partial x^{2}}(x,x_{0},x_{1})\right)^{2} dx dx_{0} dx_{1} = \frac{3 (1 + a^{2})^{2}}{64 \pi \sqrt{\pi} (1 - a^{2})^{3}} \sigma_{a}^{-7}, \\
&c^{\vt}_{0} \, =  \, \frac{1}{2} \, \kappa_{2}^{2} \iiint_{\RR^{3}} \left(\frac{\partial^{2} \mu^{\vt}}{\partial x_{0}^{2}}(x,x_{0},x_{1})\right)^{2} dx dx_{0} dx_{1}  = \frac{3}{64 \pi \sqrt{\pi} (1 - a^{2})^{3}} \sigma_{a}^{-7}, \\
&c^{\vt}_{1} \, =  \, \frac{1}{2} \, \kappa_{2}^{2} \iiint_{\RR^{3}} \left(\frac{\partial^{2} \mu^{\vt}}{\partial x_{1}^{2}}(x,x_{0},x_{1})\right)^{2} dx dx_{0} dx_{1} = \frac{3}{64 \pi \sqrt{\pi} (1 - a^{2})^{3}} \sigma_{a}^{-7}. 
\end{align*}
We have
\begin{equation*}
 \frac{c_{1}}{4 c_{2}} = b_{1} \, \sigma_{a}^{5}, \quad \frac{c_{1}}{M_{a,\sigma}} = b_{2} \, \sigma_{a}^{2}, \quad \frac{c_{2}^{\vt}}{4 \, c^{\vt}} \, =  \, b_{3} \, \sigma_{a}^{7}, \quad \frac{c_{2}^{\vt}}{4 \, c_{0}^{\vt}} = \frac{c_{2}^{\vt}}{4 \, c_{1}^{\vt}} = b_{4}  \, \sigma_{a}^{7}, \quad \frac{c_{2}^{\vt}}{M_{a,\sigma}^{\vt}} = b_{5} \; \sigma_{a}^{4}, 
\end{equation*}
where
\begin{align*}
&b_{1} = \frac{32}{3},  \quad b_{2} = \sqrt{\pi} a^{2} (1 - a^{2}) (1 + a) (1 - a)^{2} (2a^{2} - 1), \quad b_{3} = \frac{48(1 - a^{2})^{3}}{12(1 + a^{2})^{2}},  \\
&b_{4} = 4 (1-a^{2})^{3}, \quad b_{5} = \frac{3 a^{2} (1 - a^{2})^{3} (1 + a) (2 a^{2} - 1)}{4 \sqrt{\pi}}.
\end{align*}
Since for $a \in (0,1)$ the constants $b_{i}$, $i \in \{1, \ldots, 5\},$ are bounded, we can approximate  $h_{D},$ $h_{N},$ $h_{0N},$ $h_{1N}$ by:
\begin{align*}
&\widehat{h}_{D} = |\GG_{n}|^{-1/5} \, \widehat{\sigma}_{a} \, \ind_{\{2 \, \widehat{a}^{2} \leq 2 ^{-1/5}\}} \, + \, (2 \, \widehat{a}^{2})^{-n} \, \widehat{\sigma}_{a} \,  \ind_{\{2 \widehat{a}^{2} > 2 ^{-1/5}\}}, \\
&\widehat{h}_{N} = \widehat{h}_{0N} = \widehat{h}_{1N} = |\GG_{n}|^{-1/7} \, \widehat{\sigma}_{a} \, \ind_{\{2 \, \widehat{a}^{2} \leq 2 ^{-3/7}\}} \, + \, (2 \, \widehat{a}^{2})^{-n/3} \, \widehat{\sigma}_{a} \,  \ind_{\{2 \widehat{a}^{2} > 2 ^{-3/7}\}}, 
\end{align*}
where $\widehat{\sigma}_{a}$ is the estimator of the standard deviation of the measure $\mu$ and $\widehat{a}$ is the estimator of the geometric ergodic rate. Note that in practice, the estimators $\widehat{h}_{N}$, $\widehat{h}_{0N}$ and $\widehat{h}_{1N}$ differ slightly. Indeed, for $\widehat{h}_{D}$ and $\widehat{h}_{N}$, $\widehat{\sigma}_{a}$ is computed using the sample $(X_{u}, u \in \GG_{n})$, for  $\widehat{h}_{0N}$, $\widehat{\sigma}_{a}$ is computed using the sample $(X_{u0}, u \in \GG_{n})$ and for  $\widehat{h}_{1N}$, $\widehat{\sigma}_{a}$ is computed using the sample $(X_{u1}, u \in \GG_{n})$. Recall that for $i \in \TT$ and $A \subset \TT$, $iA = \{ij, j \in A\}$, where $ij$ is the concatenation of the two sequences $i, j \in \TT.$ For the geometric ergodic rate, we propose the following estimator, which is inspired from \cite{gyo-pau2016}: 
\begin{equation*}
\widehat{a} = \left(\frac{\sum_{u \in \GG_{n-m+1}}\sum_{v \in u\GG_{m-1}} (X_{u} - \overline{X})(X_{v} - \overline{X})}{\sum_{u \in \GG_{n}}(X_{u} - \overline{X})^{2}}\right)^{1/m} \quad \text{with} \quad \overline{X} = \frac{1}{|\GG_{n}|} \sum_{u \in \GG_{n}} X_{u},
\end{equation*}
where $m$ is a large enough natural integer such that $m = \bigO(n)$. The choice $m = \lfloor n/2 \rfloor + 1$ seems to be relevant. Finally, the estimator used for numerical studies is $\widetilde{\Pp}_{\GG_{n}}$ defined by:
\begin{equation*}
\widetilde{\Pp}_{\GG_{n}}(xx_{0}x_{1}) = \frac{\widetilde{\mu}^{\vartriangle}_{\A_{n}}(xx_{0}x_{1})}{\widetilde{\mu}_{\A_{n}}(x)},
\end{equation*}
with
\begin{align*}
&\widetilde{\mu}_{\GG_{n}}(x) = \frac{1}{|\GG_{n}|\widehat{h}_{D}} \sum_{u \in \GG_{n}} K_{0}\left(\frac{x - X_{u}}{\widehat{h}_{D}}\right)  \quad \text{and} \quad \\
&\widetilde{\mu}^{\vartriangle}_{\GG_{n}}(xx_{0}x_{1}) = \frac{1}{|\GG_{n}|\widehat{h}_{N}^{3}}\sum_{u \in \GG_{n}} K_{0}\left(\frac{x - X_{u}}{\widehat{h}_{N}}\right) K_{0}\left(\frac{x_{0} - X_{u0}}{\widehat{h}_{0N}}\right) K_{0}\left(\frac{x_{1} - X_{u1}}{\widehat{h}_{1N}}\right). 
\end{align*}
This method is an adaptation of the rule of thumb developed by Silverman in \cite{silverman1986density}. The novelty here is that the ergodic rate of convergence is taken into account in the estimation procedure. In the context of BMC, The main advantage of this method is that it not requires a lot of time for calculations. However, this method is a crude approximation which works for approximately ``Gaussian'' bifurcating Markov chains.

\subsection{Numerical illustrations}

In order to validate our method, we consider two cases:
\begin{enumerate}
\item[] case 1: $(a_{0}, a_{1}, b_{0}, _{1}, \sigma, \rho)= (0.7, 0.5, 0, 0, 1, 0);$
\item[] case 2: $(a_{0}, a_{1}, b_{0}, _{1}, \sigma, \rho)= (1.2, 0.7, 0, 0, 1, 0);$  
\end{enumerate}
In case 2, we allow the dynamic of the new pole to be unstable, even if the entire dynamic of the system is stable. Following the terminology of Bitseki and Delmas in \cite{BD2, BD1}, the case 2 corresponds to supercritical case.

\begin{figure}[!ht]
	\centering
	\begin{subfigure}{0.45\textwidth} 
		\includegraphics[width=\textwidth]{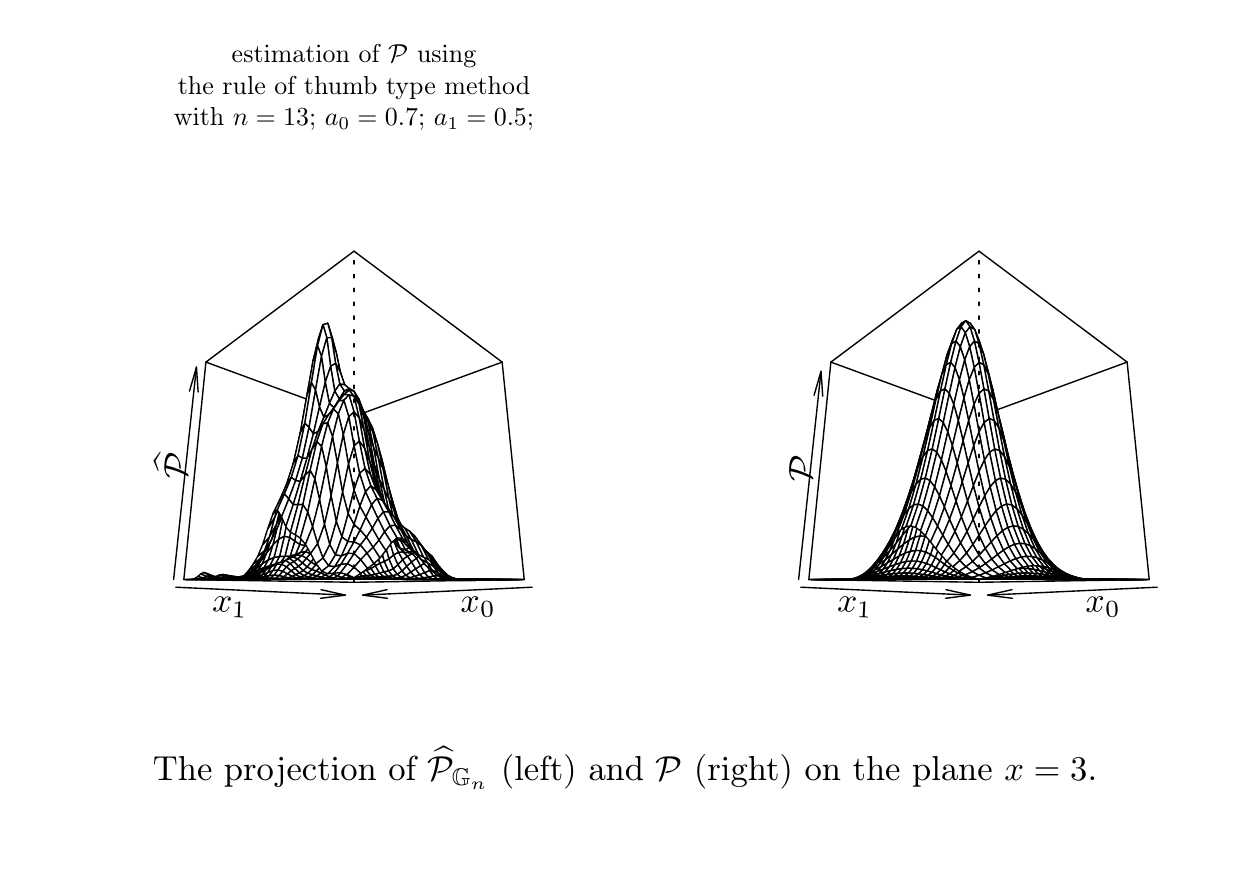}
	\end{subfigure}
	\begin{subfigure}{0.45\textwidth} 
		\includegraphics[width=\textwidth]{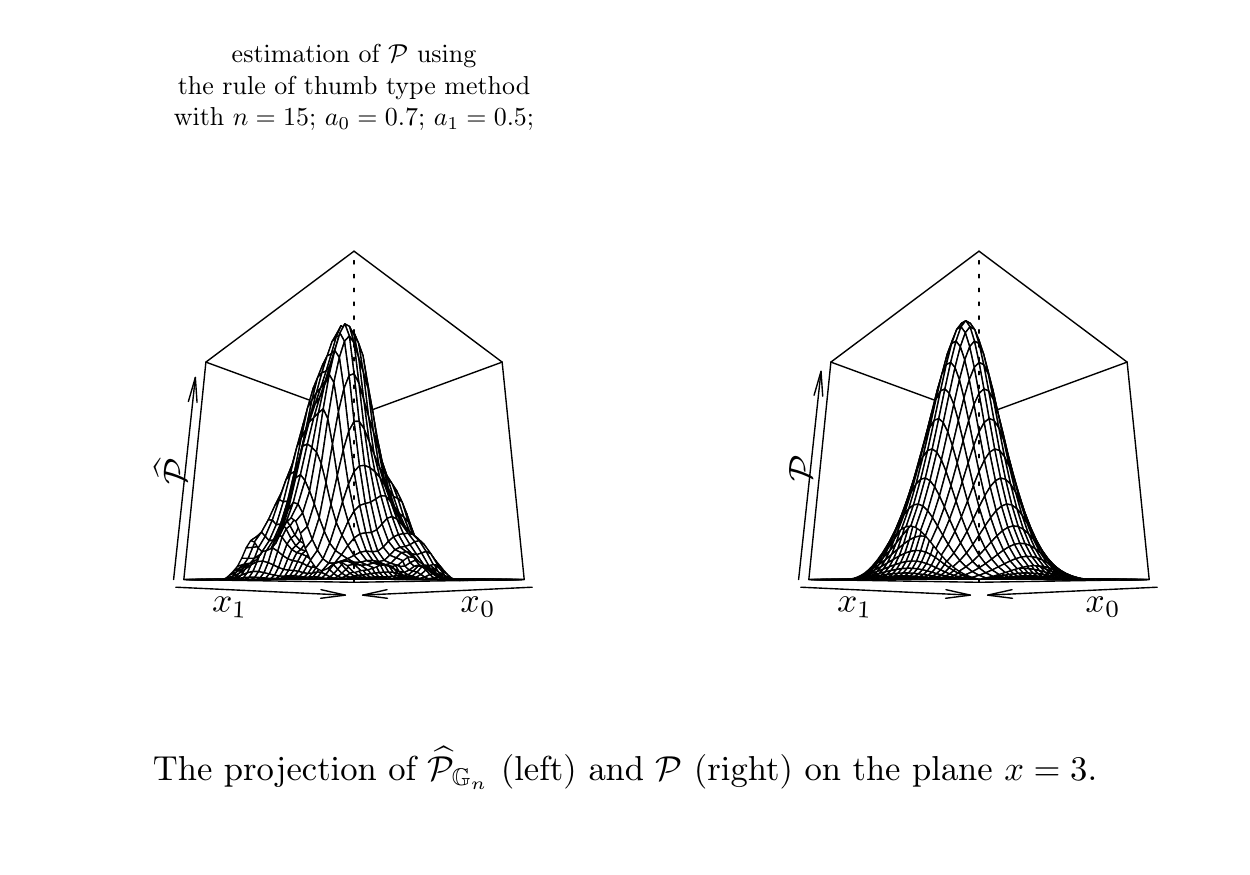}
	\end{subfigure}
        \caption{} \label{fig1}
\end{figure}

\begin{figure}[!ht]
	\centering
	\begin{subfigure}{0.45\textwidth} 
		\includegraphics[width=\textwidth]{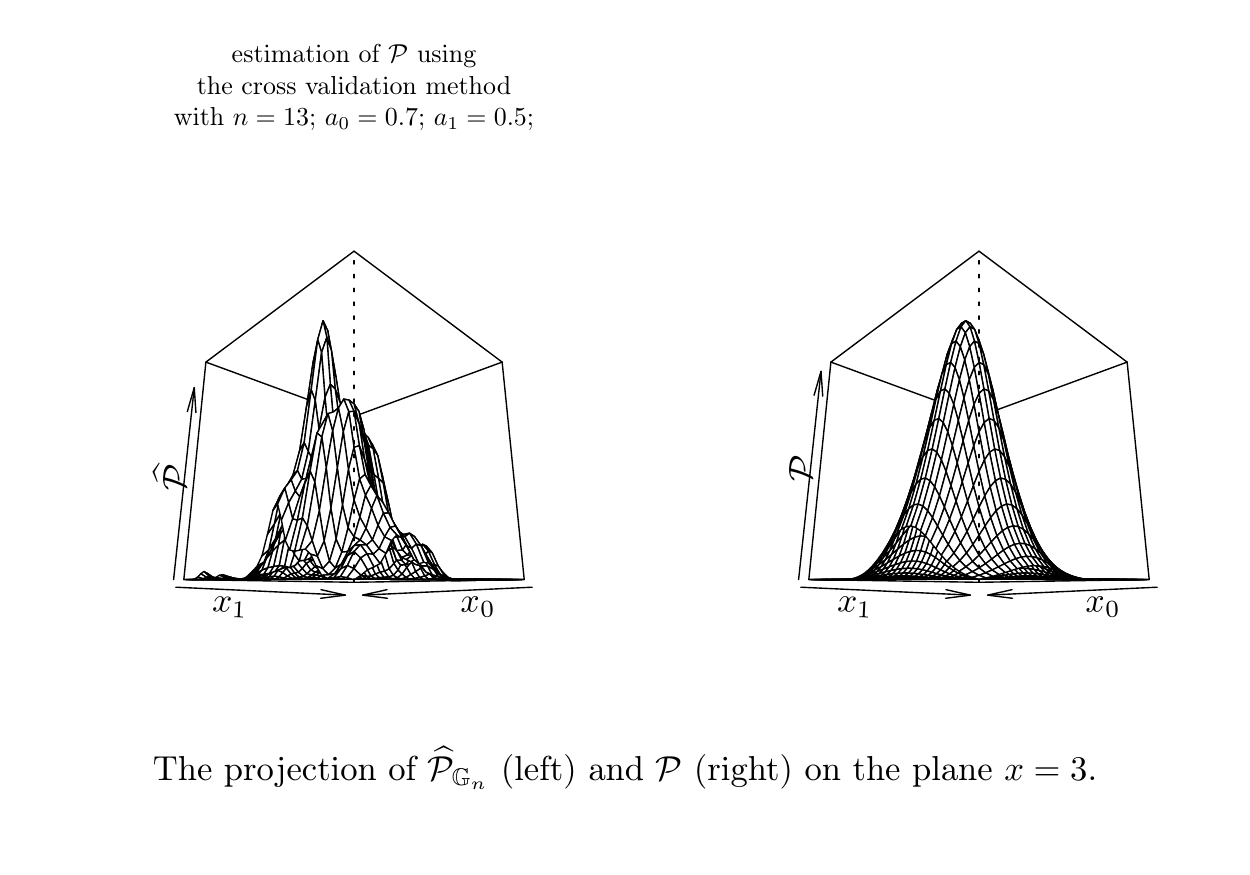}
	\end{subfigure}
	\begin{subfigure}{0.45\textwidth} 
		\includegraphics[width=\textwidth]{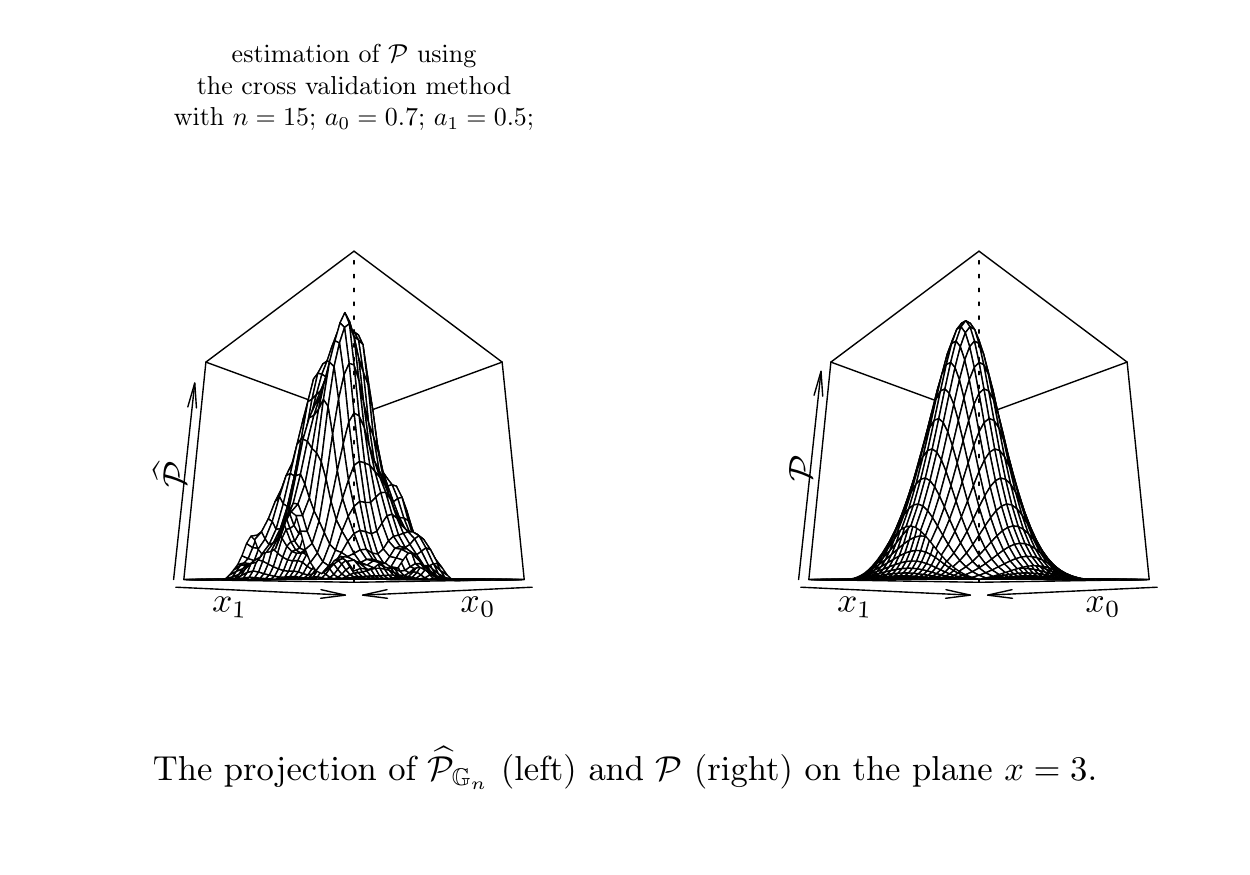}
	\end{subfigure}
        \caption{} \label{fig2}
\end{figure}

\begin{figure}[!ht]
	\centering
	\begin{subfigure}{0.45\textwidth} 
		\includegraphics[width=\textwidth]{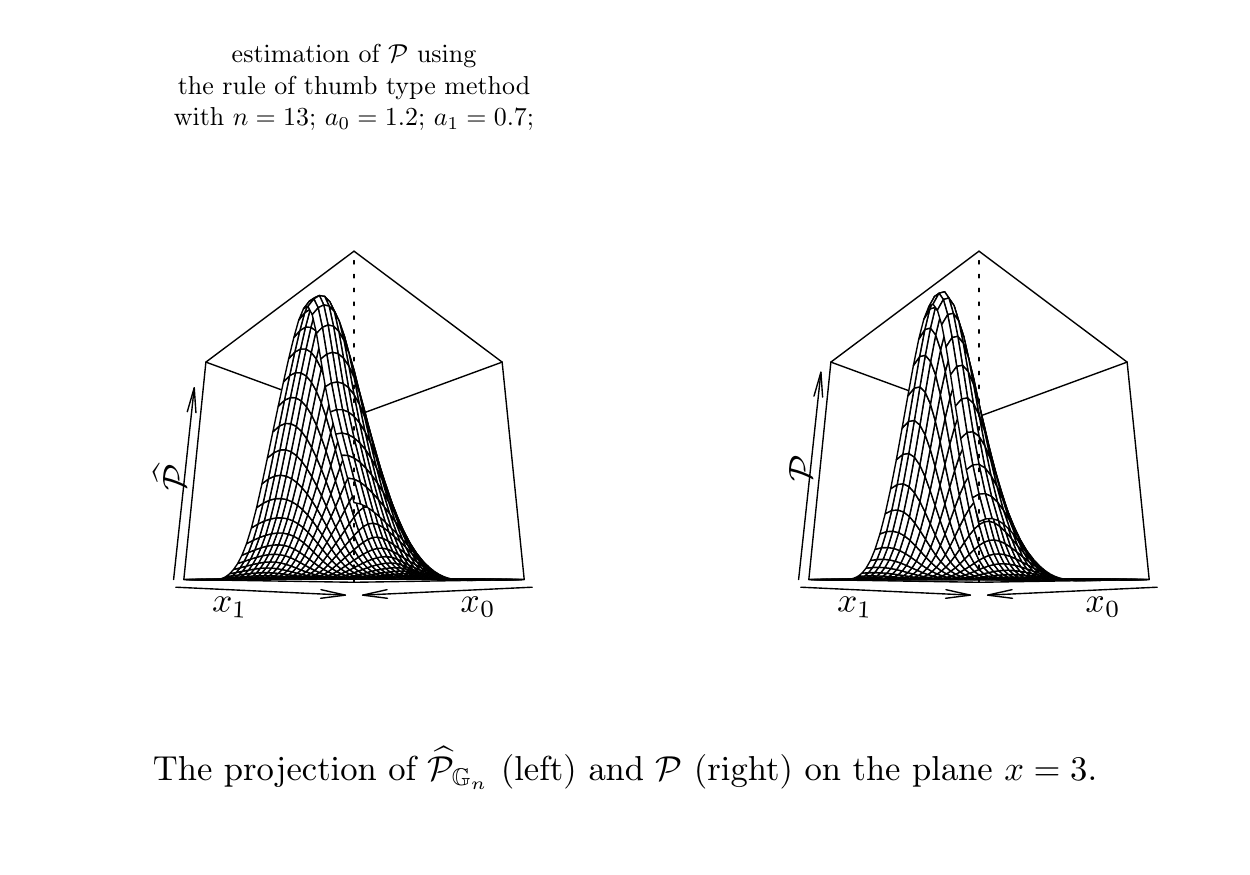}
	\end{subfigure}
	\begin{subfigure}{0.45\textwidth} 
		\includegraphics[width=\textwidth]{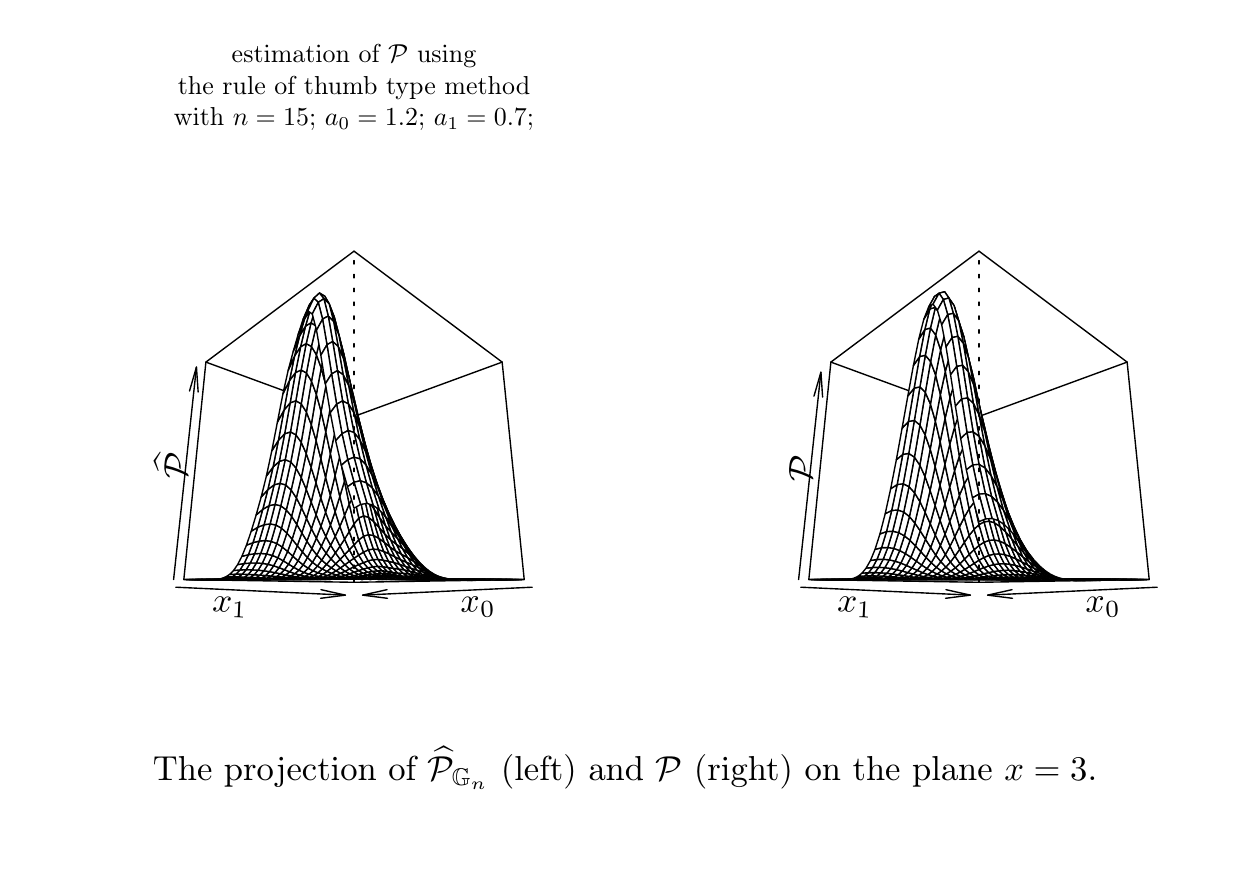}
	\end{subfigure}
        \caption{} \label{fig3}
\end{figure}

\begin{figure}[!ht]
	\centering
	\begin{subfigure}{0.45\textwidth} 
		\includegraphics[width=\textwidth]{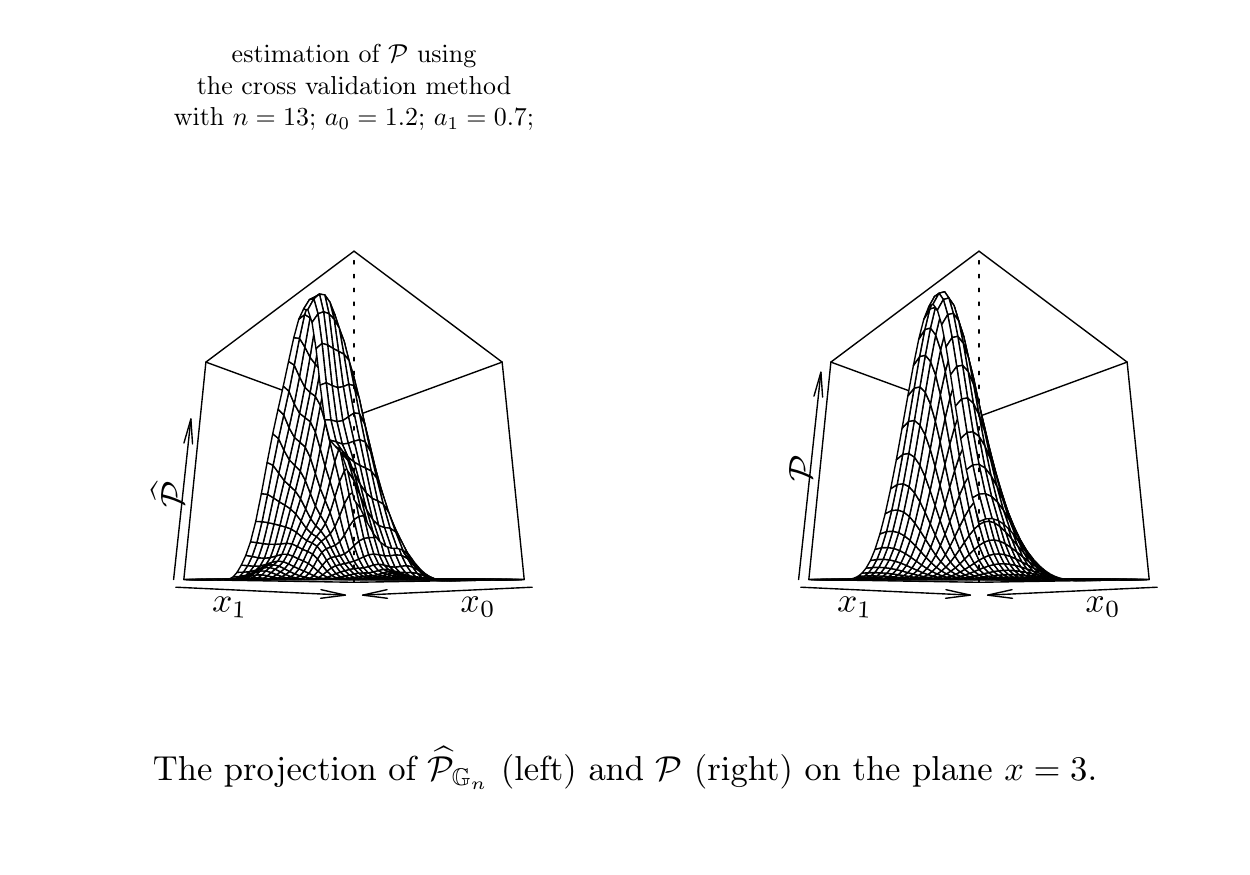}
	\end{subfigure}
	\begin{subfigure}{0.45\textwidth} 
		\includegraphics[width=\textwidth]{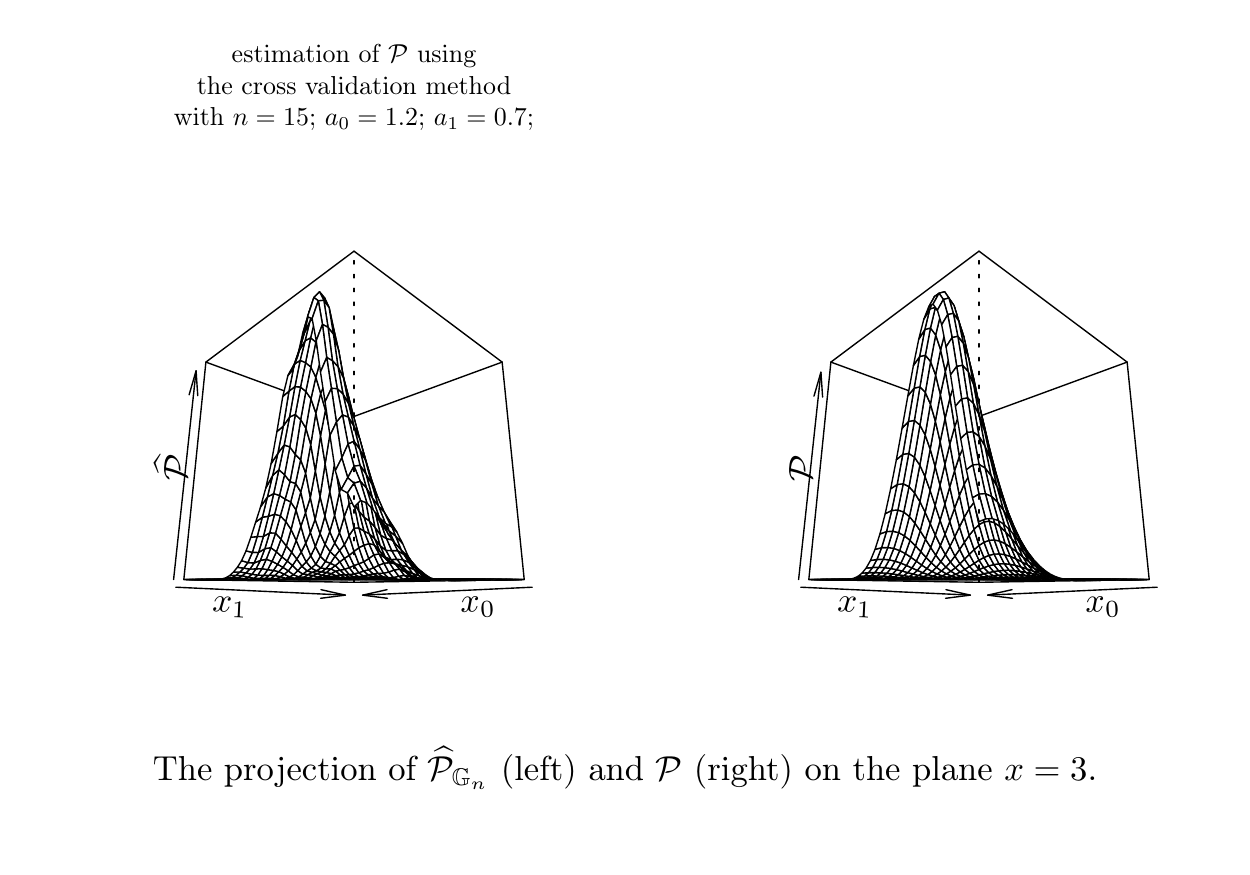}
	\end{subfigure}
        \caption{} \label{fig4}
\end{figure}

\begin{figure}[!ht]
	\centering
	\begin{subfigure}{0.45\textwidth} 
		\includegraphics[width=\textwidth]{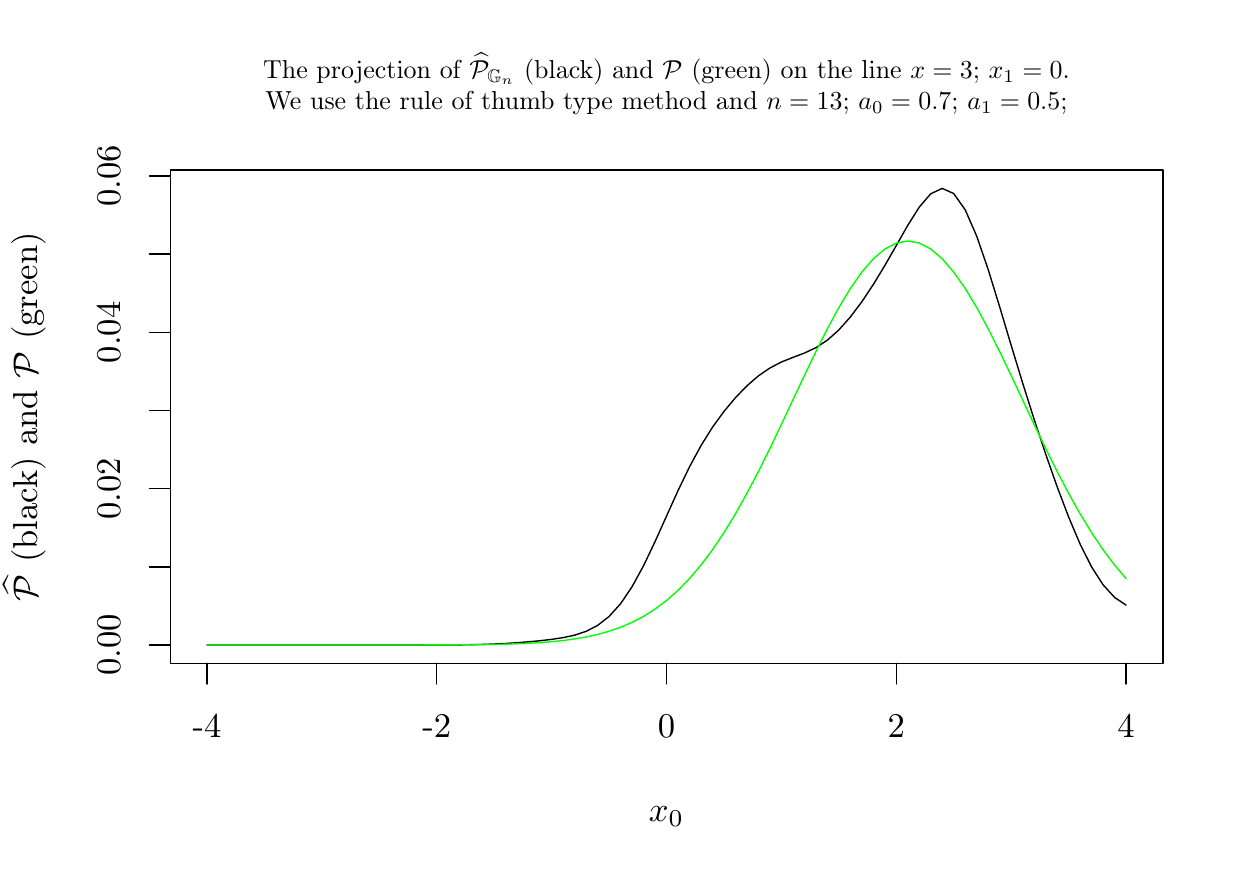}
	\end{subfigure}
	\begin{subfigure}{0.45\textwidth} 
		\includegraphics[width=\textwidth]{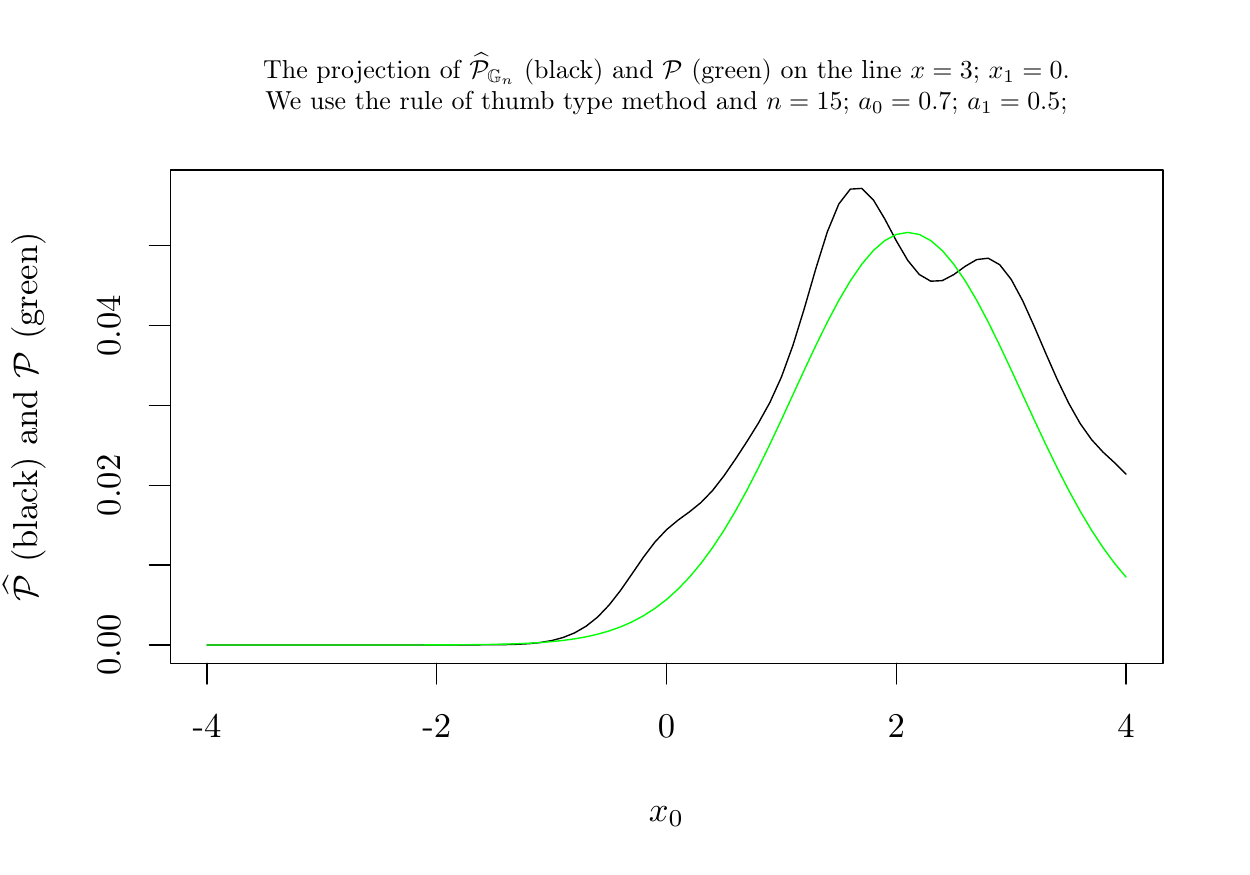}
	\end{subfigure}
        \caption{} \label{fig5}
\end{figure}

\begin{figure}[!ht]
	\centering
	\begin{subfigure}{0.45\textwidth} 
		\includegraphics[width=\textwidth]{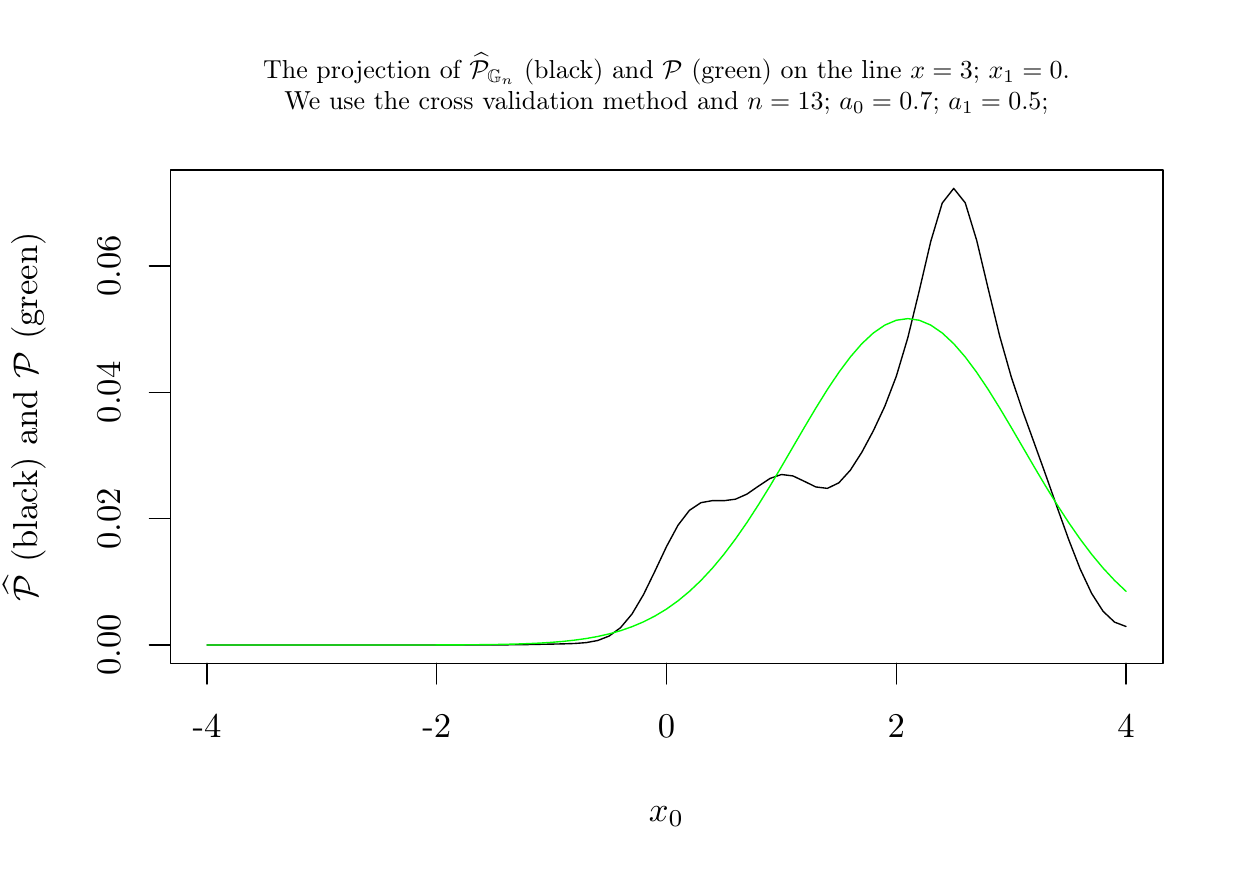}
	\end{subfigure}
	\begin{subfigure}{0.45\textwidth} 
		\includegraphics[width=\textwidth]{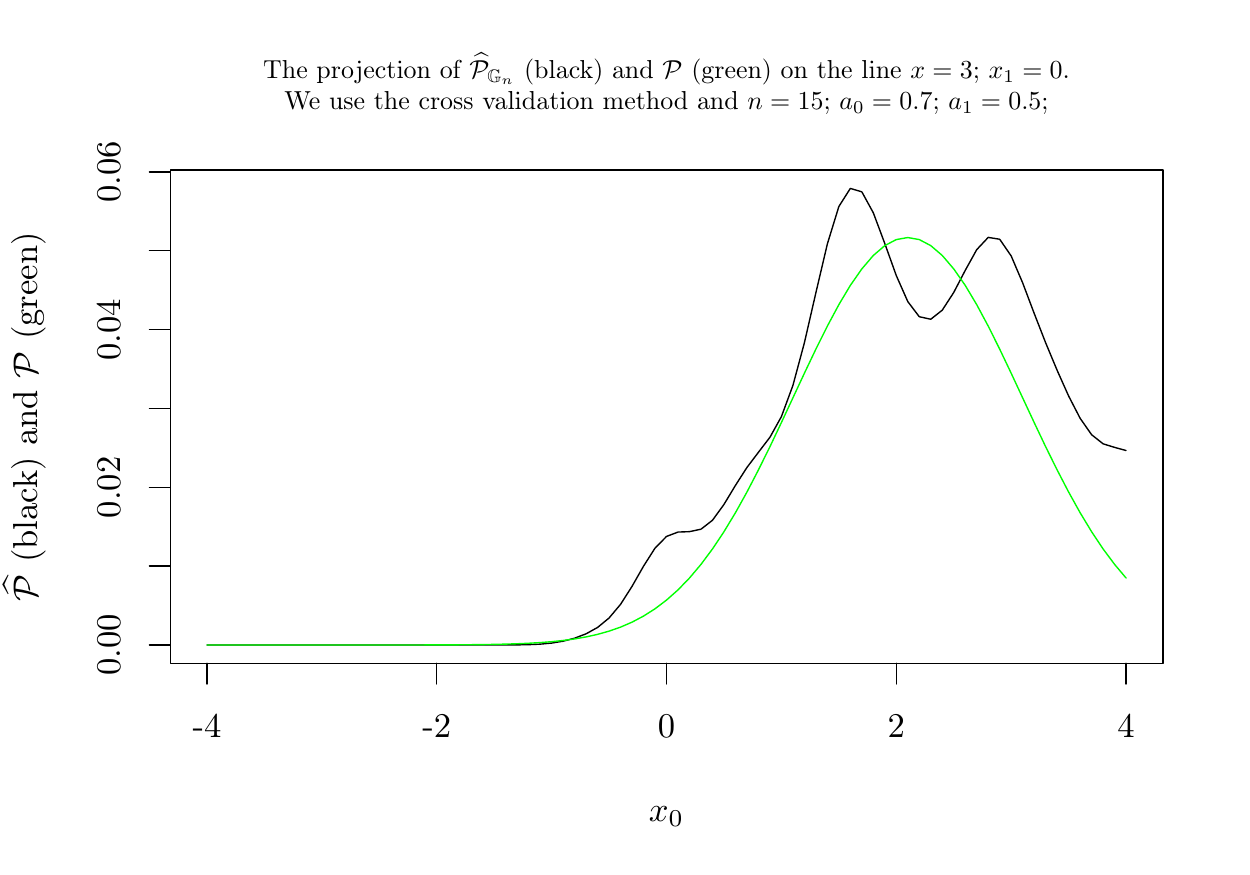}
	\end{subfigure}
        \caption{} \label{fig6}
\end{figure}

\begin{figure}[!ht]
	\centering
	\begin{subfigure}{0.45\textwidth} 
		\includegraphics[width=\textwidth]{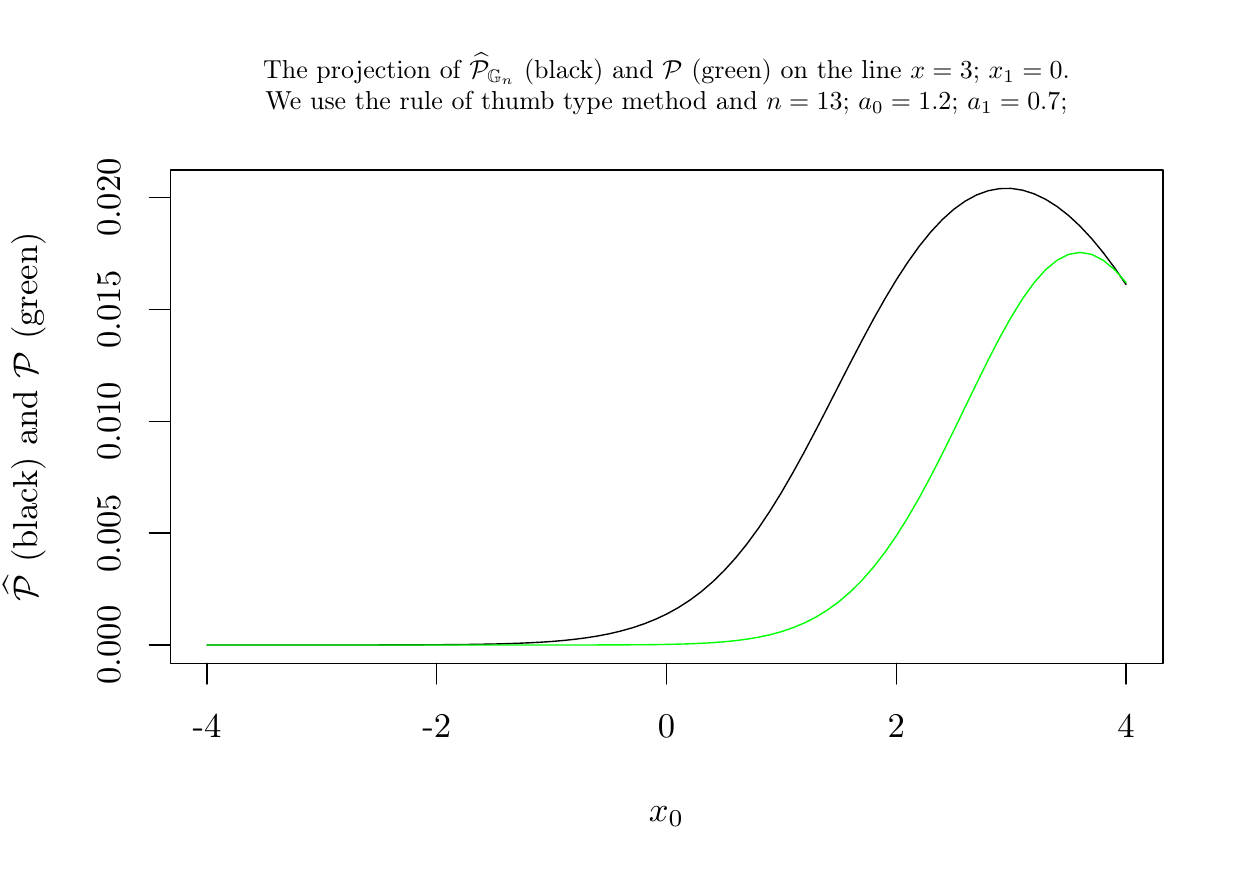}
	\end{subfigure}
	\begin{subfigure}{0.45\textwidth} 
		\includegraphics[width=\textwidth]{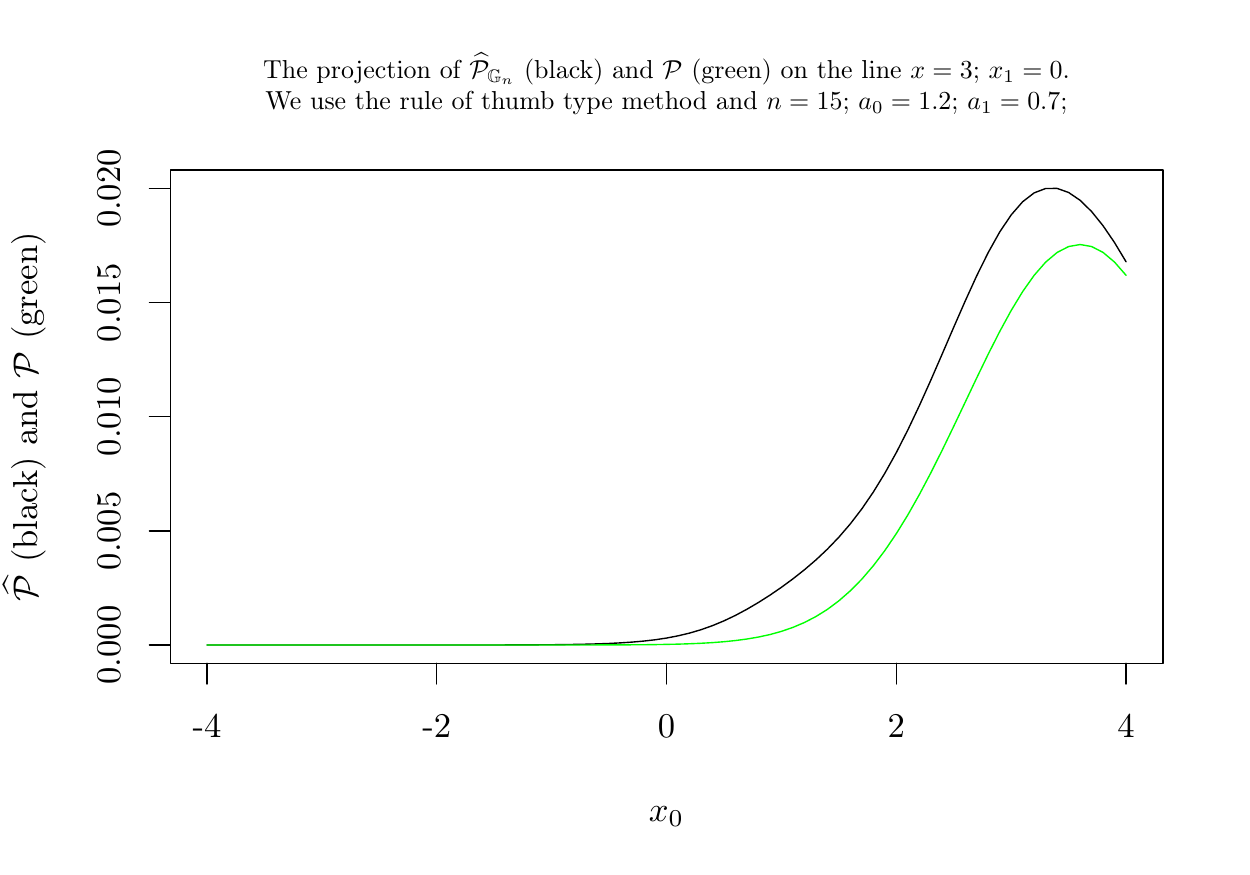}
	\end{subfigure}
        \caption{} \label{fig7}
\end{figure}

\begin{figure}[!ht]
	\centering
	\begin{subfigure}{0.45\textwidth} 
		\includegraphics[width=\textwidth]{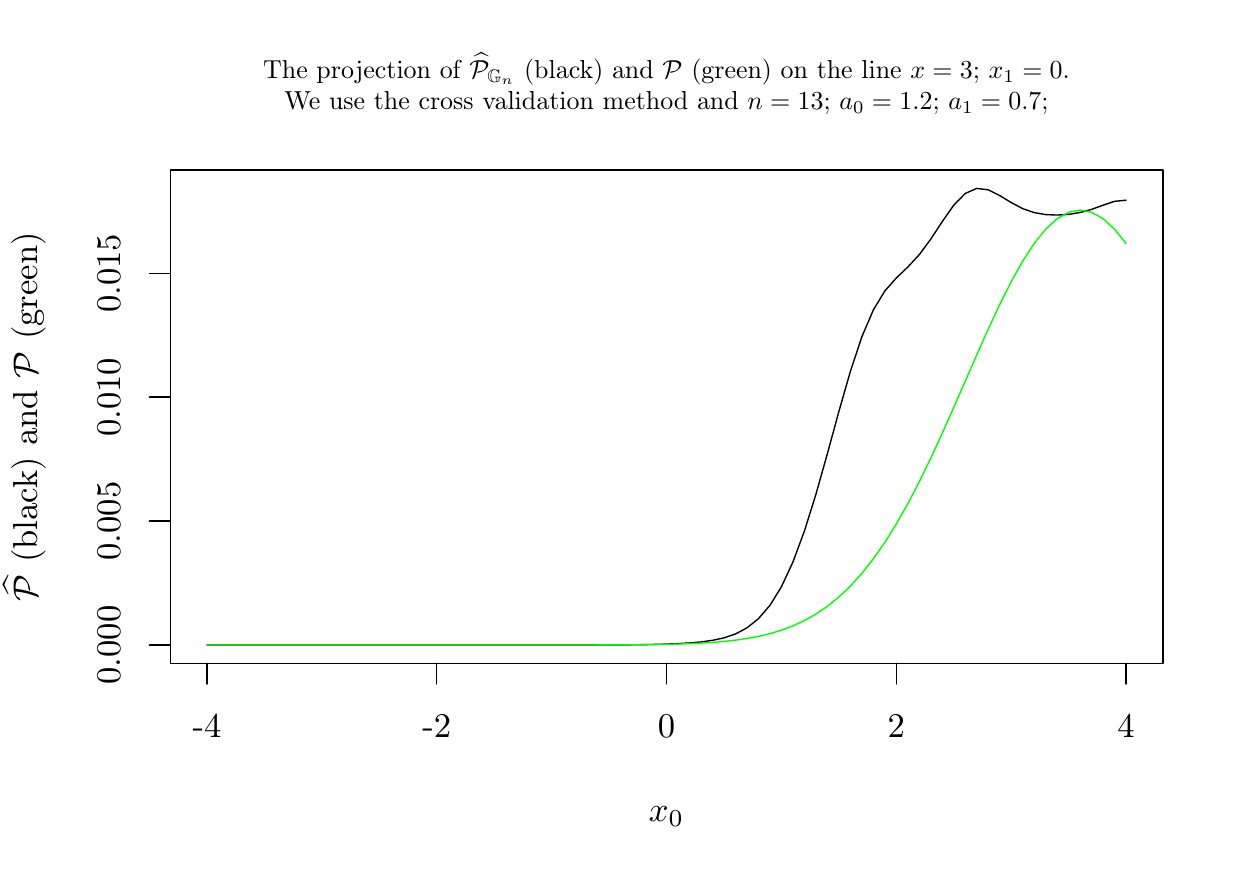}
	\end{subfigure}
	\begin{subfigure}{0.45\textwidth} 
		\includegraphics[width=\textwidth]{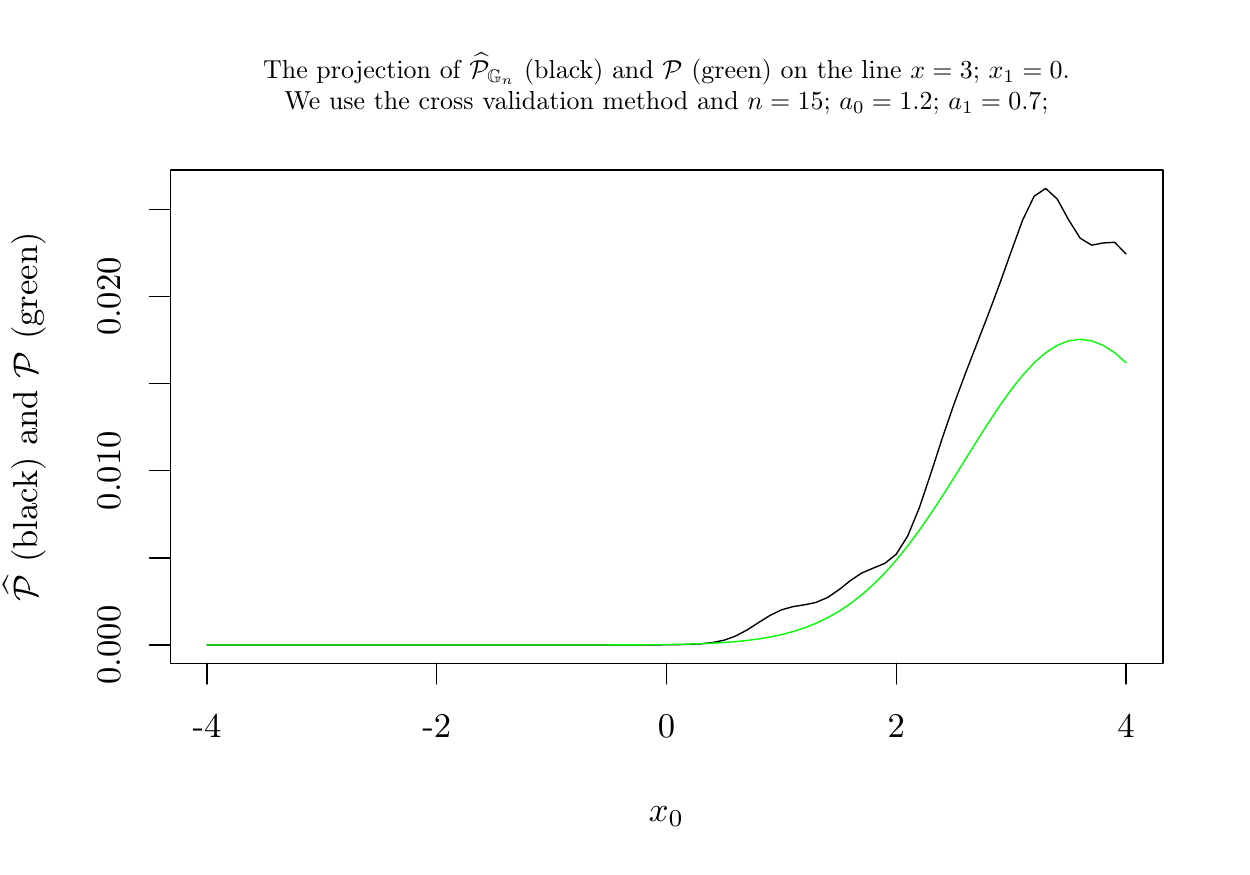}
	\end{subfigure}
        \caption{} \label{fig8}
\end{figure}

As we can see, Figure \ref{fig1}-\ref{fig8}, the two methods allow to recover the true function when the size of the data increases. Consequently, we conclude that our method is valid.

\section{Proof of Theorem \ref{thm:flx-T}}\label{proof:thm-flx-T}

We begin the proof with $\A_{n} = \TT_{n}.$ Let $(p_n, n\in \N)$ be a non-decreasing sequence of elements of $\N^*$ such that, for all $\lambda>0$:
\begin{equation*}\label{eq:def-pn}
p_n< n, \quad \lim_{n\rightarrow \infty } p_n/n=1 \quad \text{and}\quad \lim_{n\rightarrow \infty } n-p_n - \lambda \log(n)=+\infty .
\end{equation*}
When there is no ambiguity, we write $p$ for $p_n$. Recall the function $f_{n}$ defined in \eqref{eq:def-fnT}. We have the following decomposition:
\begin{equation}\label{eq:DNn0k0-T}
N_{n,\emptyset}(f_{n}) = R_{0}(n) \, + \, R_{1}(n) \, + \, \Delta_{n}(f_{n}),  
\end{equation}
where:
\begin{equation*}
R_{0}(n) = |\GG_{n}|^{-1/2} \sum_{k = 0}^{n-p-1} M_{\GG_{k}}(\tilde{f}_{n}); \; \; R_{1}(n) = \sum_{i \in \GG_{n-p}} \EE[N_{n,i}(f_{n})|\Ff_{i}]; \; \; \Delta_{n}(f_{n}) = \sum_{i \in \GG_{n-p}} \Delta_{n,i}(f_{n}),
\end{equation*}
and for all $i \in \GG_{n-p}$,
\begin{equation*}
N_{n,i}(f_{n}) = |\GG_{n}|^{-1/2} \sum_{\ell = 0}^{p} M_{i\GG_{p-\ell}}(\tilde{f}_{n}) \quad \text{and} \quad \Delta_{n,i}(f_{n}) = N_{n,i}(f_{n}) - \EE[N_{n,i}(f_{n})|\Ff_{i}].
\end{equation*}
Note that using the branching Markov property, we have, for all $i \in \GG_{n-p}$,
\begin{equation}\label{eq:EMiGp-lfn-T}
\EE[M_{i\GG_{p-\ell}}(\tilde{f}_{n})|\Ff_{i}] = \EE_{X_{i}}[M_{\GG_{p-\ell}}(\Pp \tilde{f}_{n})].
\end{equation} 
We have the following convergence.
\begin{lem}\label{lem:cge-R0-est}
Under the assumptions of Theorem \ref{thm:flx-T}, we have that $\lim_{n \rightarrow \infty} \EE[R_{0}(n)^{2}] = 0.$
\end{lem}
\begin{proof}
We have
\begin{align}
\EE[R_{0}(n)^{2}]  &= |\GG_{n}|^{-1} \, \EE[(\sum_{k = 0}^{n-p-1} \sum_{u\in \GG_{k}} \tilde{f}_{n}(X_{u}^{\vt})^{2}]& \nonumber\\ &\leq |\GG_{n}|^{-1} \, (\sum_{k = 0}^{n-p-1} \EE[(\sum_{u\in\GG_{k}}\tilde{f}_{n}(X_{u}^{\vt})^{2}]^{1/2})^{2},& \label{eq:BR0Triang}
\end{align}
where we used the Minkowski inequality for the first inequality. By developing the term in the expectation, we get
\begin{align*}
\EE[(\sum_{u\in\GG_{k}}\tilde{f}_{n}(X_{u}^{\vt})^{2}] &= \EE[\sum_{u\neq v\in\GG_{k}} \EE[\tilde{f}_{n}(X_{u}^{\vt}) \tilde{f}_{n}(X_{v}^{\vt})|X_{u},X_{v}]] + \EE[\sum_{u\in\GG_{k}} \EE[(\tilde{f}_{n})^{2}(X_{u}^{\vt})|X_{u}]]& \\ &= \EE[\sum_{u\neq v\in\GG_{k}} \cp\tilde{f}_{n}(X_{u}) \cp\tilde{f}_{n}(X_{v})] + \EE[\sum_{u\in\GG_{k}} \cp((\tilde{f}_{n})^{2})(X_{u})]& \\ & = \EE[(\sum_{u \in \GG_{k}}\cp\tilde{f}_{n}(X^{\vt}_{u}))^{2}] + \EE[\sum_{u\in\GG_{k}} (\cp((f_{n})^{2}) - (\cp f_{n})^{2})(X_{u})],& 
\end{align*}
where we used the branching Markov property for the second inequality and the fact that $\cp((\tilde{f}_{n})^{2}) - (\cp \tilde{f}_{n})^{2} = \cp((f_{n})^{2}) - (\cp f_{n})^{2}$ for the third equality. Using \eqref{eq:BR0Triang} and using the inequalities $\sqrt{a+b} \leq \sqrt{a}+\sqrt{b}$ and $(a+b)^{2} \leq 2a^{2} + 2b^{2}$, we get
\begin{multline*}\label{eq:bdRtriSub}
\EE[R_{0}(n)^{2}] \leq |\GG_{n}|^{-1} \, (\sum_{k = 0}^{n-p-1} (\EE[(\sum_{u\in\GG_{k}}\cp\tilde{f}_{n}(X^{\vt}_{u}))^{2}]^{1/2} + \EE[\sum_{u\in\GG_{k}} (\cp((f_{n})^{2}) - (\cp f_{n})^{2})(X^{\vt}_{u})]^{1/2}))^{2} \\ 
\leq 2 |\GG_{n}|^{-1}( (\sum_{k = 0}^{n-p-1} \EE[M_{\GG_{k}}(\cp(\tilde{f}_{n}))^{2}]^{1/2})^{2} + (\sum_{k = 0}^{n-p-1}\EE[M_{\GG_{k}}(\cp((f_{n})^{2}) - (\cp f_{n})^{2})]^{1/2})^{2}).
\end{multline*}
Note that from Lemma \ref{lem:useful-upb}, we have $\|\Pp(f_{n})\|_{\infty} \leq C h_{n}^{d/2}$ and $\|\Pp(f_{n}^{2})\|_{\infty} \leq C$. Recall $\Pp(\tilde{f}_{n}) = \Pp(f_{n}) - \langle \mu,\Pp(f_{n}) \rangle.$ Then, using \eqref{eq:Q2}, \eqref{eq:geom-ergB} and Lemma \ref{lem:useful-upb}, we get
\begin{equation*}
\EE[M_{\GG_{k}}(\Pp\tilde{f}_{n})^{2}] \leq C \, h_{n}^{d} \quad \text{if $k \in \{0,1\}$} 
\end{equation*}
and for all $k \geq 2,$
\begin{equation*}
\EE[M_{\GG_{k}}(\Pp\tilde{f}_{n})^{2}] \leq \begin{cases} C \, h_{n}^{2d} \, 2^{k} & \text{if $2\alpha^{2} \leq 1$} \\ C \, 2^{k} \,  (h_{n}^{2d} \, + \, (2\alpha^{2})^{k} \, h_{n}^{3d}) & \text{if $2\alpha^{2} > 1.$} \end{cases}
\end{equation*}
It follows for the two last inequalities that
\begin{equation*}
|\GG_{n}|^{-1} \, (\sum_{k = 0}^{n-p-1} \EE[M_{\GG_{k}}(\cp(\tilde{f}_{n}))^{2}]^{1/2})^{2} \leq  C \, 2^{-n} \, h_{n}^{3d} \, + \, C \, 2^{-p} \,  h_{n}^{2d} \, + \, 2^{-p} \, (2\alpha^{2})^{n-p} \, h_{n}^{3d} \,\ind_{\{2\alpha^{2} > 1\}}.
\end{equation*} 
Using \eqref{eq:2alpha2h3d}, Assumption \ref{hyp:2alpha2h3d} and since $\lim_{n \rightarrow \infty} p_{n} = \infty$, it follows that
\begin{equation*}
\lim_{n \rightarrow \infty} |\GG_{n}|^{-1}(\sum_{k = 0}^{n-p-1} \EE[M_{\GG_{k}}(\cp(\tilde{f}_{n}))^{2}]^{1/2})^{2} = 0.
\end{equation*}
Next, using Lemma \ref{lem:useful-upb}, we get
$\EE[M_{\GG_{k}}(\cp((f_{n})^{2}) - (\cp f_{n})^{2})] \leq C \, 2^{k}.$
This implies that
\begin{equation*}
\lim_{n \rightarrow \infty}|\GG_{n}|^{-1}(\sum_{k = 0}^{n-p-1}\EE[M_{\GG_{k}}(\cp((f_{n})^{2}) - (\cp f_{n})^{2})]^{1/2})^{2} \, \leq \, C \, \lim_{n \rightarrow \infty} 2^{-p} \, = \, 0
\end{equation*}
and this ends the proof.
\end{proof}

Next, we have the following convergence.

\begin{lem}\label{lem:cge-R1-est}
Under the assumptions of Theorem \ref{thm:flx-T}, we have that $\lim_{n \rightarrow \infty} \EE[R_{1}(n)^{2}] = 0.$
\end{lem}

\begin{proof}
Using \eqref{eq:EMiGp-lfn-T}, we get
\begin{equation*}
R_{1}(n) = \sum_{k = 0}^{p} R_{1}(k,n),
\end{equation*}
with
\begin{equation*}
R_{1}(k,n) = |\G_n|^{-1/2} \,  |\G_{p-k}|\, M_{\G_{n-p}}
(\cq^{p-k} \Pp \tilde{f}_n).
\end{equation*}
It follows that
\begin{equation}\label{eq:IR1nL2sub}
\E\left[R_1(n)^2\right]^{1/2} \leq  \sum_{k=0}^{p}  \left(\E\left[R_{1}(k, n)^2\right] \right)^{1/2}.
\end{equation}
Following the proof of Lemma 4.2 in \cite{BD1} and using \eqref{eq:geom-ergB} and Lemma \ref{lem:useful-upb}, we find that
\begin{equation*}
\E\left[R_{1}(k, n)^2\right] \leq C \, 2^{-p} \, h_{n}^{2d} \, \ind_{\{k=p\}} + \begin{cases} C \, h_{n}^{3d} \, 2^{-p} \, (2\alpha)^{2(p-k)} & \text{if $2\alpha^{2} < 1$} \\ C \, (n-p) \, h_{n}^{3d} \, 2^{-k} & \text{if $2\alpha^{2} = 1$} \\ C \, 2^{-p} \, (2\alpha^{2})^{n-p} \, h_{n}^{3d} (2\alpha)^{2(p-k)} & \text{if $2\alpha^{2} > 1.$} \end{cases} 
\end{equation*}
From \eqref{eq:IR1nL2sub}, this implies that
\begin{equation*}
\E\left[R_1(n)^2\right]^{1/2} \leq C \, 2^{-p/2} \, h_{n}^{d} \, + \begin{cases} C \, h_{n}^{3d/2} \, \sum_{k=0}^{p} 2^{-k/2} \, (2\alpha^{2})^{(p-k)/2} & \text{if $2\alpha^{2} < 1$} \\ C \, (n-p)^{1/2} \, h_{n}^{3d/2} & \text{if $2\alpha^{2} = 1$} \\ C \, (2\alpha^{2})^{n/2} \, h_{n}^{3d/2} & \text{if $2\alpha^{2} > 1.$} \end{cases} 
\end{equation*}
From the latter inequality and using Assumption \ref{hyp:2alpha2h3d}, we deduce that $\lim_{n \rightarrow \infty} \EE[R_{1}(n)^{2}] = 0.$
\end{proof}

We now study the bracket 
\begin{equation*}\label{eq:bracket-tri}
V(n) = \sum_{i\in\GG_{n-p}}\EE[\Delta_{n,i}(f_{n})^{2}|\Ff_{i}].
\end{equation*} 
Note that for $i\in\GG_{n-p}$, we have
\begin{equation*}
\EE[\Delta_{n,i}(f_{n})^{2}|\Ff_{i}] = |\GG_{n}|^{-1} \EE[(\sum_{k=0}^{p} M_{i\GG_{p - k}}(\tilde{f}_{n}))^{2}|\Ff_{i}] - |\GG_{n}|^{-1}(\EE[\sum_{k = 0}^{p} M_{i\GG_{p-k}}(\tilde{f}_{n})|\Ff_{i}])^{2}.
\end{equation*}
Using the branching Markov chain property, this implies that 
\begin{equation}\label{eq:dVtriang}
V(n) = V_{1}(n) + V_{2}(n) - R_{2}(n),
\end{equation} 
with
\begin{align}
&V_{1}(n) = |\GG_{n}|^{-1}\sum_{i \in \GG_{n-p}}\sum_{k = 0}^{p} \EE_{X_{i}}[M_{\GG_{p - k}}(\tilde{f}_{n})^{2}],& \nonumber \\ &V_{2}(n) = 2 |\GG_{n}|^{-1}\sum_{i\in\GG_{n-p}}\sum_{0\leq k < \ell \leq p} \EE_{X_{i}}[M_{\GG_{p-\ell}}(\tilde{f}_{n})M_{\GG_{p-k}}(\tilde{f}_{n})],& \label{eq:V2n-est-tri} \\ &R_{2}(n) = |\GG_{n}|^{-1} \sum_{i \in \GG_{n - p}} (\sum_{k = 0}^{p} 2^{p-k} \Qq^{p-k} \cp\tilde{f}_{n}(X_{u}))^{2}.& \nonumber
\end{align}

We have the following result.

\begin{lem}\label{lem:cvR2est}
Under the assumptions of Theorem \ref{thm:flx-T}, we have the following convergence:
\begin{equation*}
\lim_{n\rightarrow \infty} \EE[R_{2}(n)] = 0.
\end{equation*}
\end{lem}
\begin{proof}
We have using \eqref{eq:Q1}, \eqref{eq:geom-ergB} and Lemma \ref{lem:useful-upb}:
\begin{align}
\nonumber
\E\left[R_2(n)\right]  &= |\G_n|^{-1} \,  |\G_{n-p}|\, \langle \nu,\cq^{n-p}\left(\Big(\sum_{k=0}^{p}  |\G_{p-k}|\, \cq^{p-k} \Pp \tilde f_n\Big)^2\right) \rangle \\
&\leq C \, 2^{-p} \langle \nu, \Qq^{n-p}((\Pp \tilde{f}_{n})^{2})\rangle + C \, 2^{-p} \langle \nu, \Qq^{n-p}((\sum_{k=0}^{p-1} 2^{p-k} \Qq^{p-k-1}(\Qq\Pp\tilde{f}_{n}))^{2}) \rangle \nonumber \\ 
& \leq \, C\, 2^{-p} \, h_{n}^{2d} \, + \, C \, 2^{-p} \, h_{n}^{3d} \, a_{n} \nonumber, 
\end{align}
where the sequence $(a_{n}, n\geq 1)$ is defined by
\begin{equation*}
a_{n} = \begin{cases} 1 & \text{if $2\alpha < 1$} \\ p^{2} & \text{if $2\alpha = 1$} \\ (2\alpha)^{2p} & \text{if $2\alpha > 1$}.  \end{cases}
\end{equation*}
Using Assumption \ref{hyp:2alpha2h3d}, and in particular Remark \ref{rem:2alpha2h3d}, it follows that $\lim_{n \rightarrow \infty} \EE[R_{2}(n)] = 0.$
\end{proof}

Next, we have the following result.

\begin{lem}\label{lem:cvV2est}
Under the assumptions of Theorem \ref{thm:flx-T}, we have the following convergence:
\begin{equation*}
\lim_{n\rightarrow \infty} \EE[V_{2}(n)^{2}] = 0.
\end{equation*}
\end{lem}

\begin{proof}
Let $0\leq k < \ell \leq p$ and $i \in \GG_{n-p}$. Conditioning two times, first by $\Hh_{i,p - k}$ and next by $\Hh_{i,p - \ell + 1}$, and using the branching Markov property, we get
\begin{equation*}
\EE_{X_{i}}[M_{\GG_{p-\ell}}(\tilde{f}_{n})M_{\GG_{p-k}}(\tilde{f}_{n})] = 2^{\ell - k - 1} \EE_{X_{i}}[M_{\GG_{p - \ell}}(\tilde{f}_{n}) M_{\GG_{p - \ell}}(g_{k,\ell,n})],
\end{equation*} 
where we set $g_{k,\ell,n} = \Qq^{\ell - k - 1} \Pp \tilde{f}_{n} \oplus \Qq^{\ell - k - 1} \Pp \tilde{f}_{n}$. Next, conditioning by $\Hh_{i,p - \ell}$ and using the branching Markov property, we get
\begin{align*}
\EE_{X_{i}}[M_{\GG_{p - \ell}}(\tilde{f}_{n}) M_{\GG_{p - \ell}}(g_{k,\ell,n})] & = \EE_{X_{i}}[M_{\GG_{p-\ell}}(\Pp(\tilde{f}_{n} g_{k,\ell,n}) - \Pp \tilde{f}_{n} \Pp g_{k,\ell,n})] \\ & \hspace{2cm} + \EE_{X_{i}}[M_{\GG_{p-\ell}}(\Pp \tilde{f}_{n}) M_{\GG_{p - \ell}}(\Pp g_{k,\ell,n})].
\end{align*}
From the foregoing and using \eqref{eq:Q1}, \eqref{eq:Q2-bis} and \eqref{eq:V2n-est-tri}, it follows that:
\begin{equation*}
V_{2}(n) = V_{5}(n) + V_{6}(n),
\end{equation*}
where
\begin{equation*}
V_{5}(n) = |\GG_{n-p}|^{-1} M_{\GG_{n-p}}(H_{5,n}) \quad \text{and} \quad V_{6}(n) = |\GG_{n-p}|^{-1} M_{\GG_{n-p}}(H_{6,n}),
\end{equation*}
with
\begin{align*}
H_{5,n} = \sum_{0 \leq k < \ell} &2^{-k} \Qq^{p - \ell}(\Pp(\tilde{f}_{n} \, g_{k,\ell,n}))\ind_{\{\ell \leq p\}} \quad \text{and} \\ & H_{6,n} =  \sum_{\substack{0\leq k < \ell \\ r \geq 0}} 2^{-k + r} \Qq^{p - \ell - r - 1} \Pp(\Qq^{r} \Pp \tilde{f}_{n} \sot \Qq^{r}\Pp g_{k,\ell,n}) \ind_{\{r + \ell < p\}}.
\end{align*}
First, we treat the term $V_{6}(n).$ Note that we have 
\begin{equation*}
\Qq^{r} \Pp g_{k,\ell,n} = 2 \, \Qq^{r + \ell - k} \Pp \tilde{f}_{n}.
\end{equation*}
We set
\begin{align*}
h_{k,\ell,r}^{(n)} = 2^{r - k + 1} \Qq^{p - 1 - (r + \ell)} &\Pp(\Qq^{r} \Pp \tilde{f}_{n} \sot \Qq^{r+\ell-k} \Pp \tilde{f}_{n}) \quad \text{and} \\ & \hspace{1cm} h_{k,\ell,r} = 2^{r - k + 1} \langle \mu, \Pp(\Qq^{r} \Pp \tilde{f}_{n} \sot \Qq^{r+\ell-k} \Pp \tilde{f}_{n})\rangle.
\end{align*}
We consider the following sums:
\begin{equation*}
H_{6}^{[n]} = \sum_{\substack{0 \leq k < \ell \\ r \geq 0}} h_{k,\ell,r} \ind_{\{r + \ell < p\}} \quad \text{and} \quad A_{6,n} = H_{6,n} - H_{6}^{[n]} =  \sum_{\substack{0 \leq k < \ell \\ r \geq 0}} (h_{k,\ell,r}^{(n)} - h_{k,\ell,r})  \ind_{\{r + \ell < p\}}. 
\end{equation*}
Using Lemma \ref{lem:useful-upb}, we have for all $0 \leq k < \ell$ and $r \geq 0$:
\begin{align}
|\Pp(\Qq^{r} \Pp \tilde{f}_{n} \sot \Qq^{r+\ell-k} \Pp \tilde{f}_{n})| & \leq \|\Qq^{r+\ell-k} \Pp \tilde{f}_{n})\|_{\infty} \, \|\Pp(\Qq^{r} \Pp \tilde{f}_{n} \sot \mathbf{1})\|_{\infty} \nonumber \\ 
& \leq \|\Qq^{r+\ell-k} \Pp \tilde{f}_{n})\|_{\infty} \, \|\Qq^{r + 1} \Pp \tilde{f}_{n}\|_{\infty}  \leq C \, h_{n}^{3d}. \label{eq:I-PQA6n1}
\end{align}
Moreover, using \eqref{eq:geom-ergB} and Lemma \ref{lem:useful-upb}, we have for all $0 \leq k < \ell$ and $r \geq 1$:
\begin{align}
|\Pp(\Qq^{r} \Pp \tilde{f}_{n} \sot \Qq^{r+\ell-k} \Pp \tilde{f}_{n})| &\leq C \, \alpha^{2r + \ell - k} \, \|\Qq(\Pp f_{n})\|_{\infty}^{2} \, \Pp(V \sot V) \nonumber \\ 
& \leq \, C \,  \alpha^{2r + \ell - k} \, h_{n}^{3d} \, \Pp(V \sot V). \label{eq:I-PQA6n2}
\end{align}
Distinguishing the cases $r = 0$ and $r \geq 1$ and using \eqref{eq:I-PQA6n1}, \eqref{eq:I-PQA6n2}, $(iv)$ of Assumption \ref{hyp:F} and \eqref{eq:erg-bd}, we get, for some $g_{1}, g \in F$,
\begin{align}
|H_{6,n} - H_{6}^{[n]}| & = \sum_{\substack{0 \leq k < \ell}} |h_{k,\ell,0}^{(n)} - h_{k,\ell,0}| \, \ind_{\{\ell < p\}} \, + \, \sum_{\substack{0 \leq k < \ell \\ r \geq 1}} |h_{k,\ell,r}^{(n)} - h_{k,\ell,r}| \, \ind_{\{r + \ell < p\}} \nonumber \\
& \leq \, C \, \sum_{\substack{0 \leq k < \ell}} 2^{-k} \alpha^{p-\ell-1} \, \|\Pp(\Pp \tilde{f}_{n} \sot \Qq^{\ell-k} \Pp \tilde{f}_{n})\|_{\infty}  \, V \nonumber \\
& \hspace{1cm} + \, C \, h_{n}^{3d} \, \sum_{\substack{0 \leq k < \ell \\ r \geq 1}} 2^{r-k} \, \alpha^{2r + \ell - k} \, \Qq^{p-\ell-r-1}\Pp(V \sot V) \, \ind_{\{r + \ell < p\}} \nonumber \\
& \leq C \, h_{n}^{3d} \, (V \, + \,  (\sum_{\substack{0 \leq k < \ell \\ r \geq 1}} 2^{r-k} \, \alpha^{2r + \ell - k} \, \ind_{\{r + \ell < p\}}) \, g_{1}) \, \leq C \, h_{n}^{3d} \, a_{n} \, g, \label{eq:I-absA6n}
\end{align} 
where
\begin{equation*}
a_{n} = \begin{cases} 1 & \text{if $2\alpha^{2} < 1$} \\ p & \text{if $2\alpha = 1$} \\ (2\alpha^{2})^{p} & \text{if $2\alpha^{2} > 1.$} \end{cases}
\end{equation*}
Using \eqref{eq:I-absA6n}, we find that
\begin{equation*}
|V_{6}(n) - H_{6}^{[n]}| \leq |\GG_{n-p}|^{-1} \, M_{\GG_{n-p}}(|H_{6}(n) - H_{6}^{[n]}|) \leq  C \, a_{n} \, h_{n}^{3d} \, |\GG_{n-p}|^{-1} \, M_{\GG_{n-p}}(g).
\end{equation*}
Using \eqref{eq:lfgn-G}, \eqref{eq:2alpha2h3d} and that $g \in L^{1}(\mu)$, we get
\begin{equation*}
\lim_{n \rightarrow +\infty} |V_{6}(n) - H_{6}^{[n]}| = 0 \quad \text{a.s. and in $L^{2}.$}
\end{equation*} 
Next, as for \eqref{eq:I-absA6n}, using \eqref{eq:I-PQA6n1}, \eqref{eq:I-PQA6n2} and that $F \subset L^{1}(\mu)$, we find that $|H_{6}^{[n]}| \leq C \, a_{n} \, h_{n}^{3d}.$ Using \eqref{eq:2alpha2h3d}, we get $\lim_{n \rightarrow + \infty} H_{6}^{[n]} = 0.$ Now, since we can write $V_{6}(n) = (V_{6}(n) - H_{6}^{[n]}) + H_{6}^{[n]},$ we conclude that $\lim_{n \rightarrow \infty} \EE[V_{6}(n)^{2}] = 0.$ 

\medskip

Next, we treat the term $V_{5}(n)$. We have $V_{5}(n) = (V_{5}(n) - H_{5}^{[n]}) + H_{5}^{[n]}$, where
\begin{equation*} 
H_{5}^{[n]} = \sum_{0 \leq k < \ell} 2^{-k} \langle \mu, \Pp(\tilde{f}_{n} \, g_{k,\ell,n})\rangle \ind_{\{\ell \leq p\}}.
\end{equation*}
Using Lemma \ref{lem:useful-upb}, we get, for all $0 \leq k < \ell$,
\begin{equation*}
|\Pp(\tilde{f}_{n}g_{k, \ell,n}) - \langle \mu, \Pp(\tilde{f}_{n}g_{k, \ell,n})\rangle| \leq \, C \, \|\Qq^{\ell - k - 1} \Pp\tilde{f}_{n}\|_{\infty} \, \|\Pp\tilde{f}_{n}\|_{\infty} \, \leq \, C \, h_{n}^{d} \ind_{\{k = \ell - 1\}} + C\, h_{n}^{2d} \, \ind_{\{k \leq \ell - 2\}}.
\end{equation*}
Using the latter inequality and distinguishing the cases $\ell = p$ and $\ell \leq p - 1$, we find that
\begin{equation*}
|V_{5}(n) - H_{5}^{[n]}| \, \leq C \, (2^{-p} h_{n}^{d} \, + \, h_{n}^{2d} \, + \, p \, h_{n}^{2d}).
\end{equation*}
This implies that $\lim_{n \rightarrow +\infty} |V_{5}(n) - H_{5}^{[n]}| = 0$ a.s. and in $L^{2}.$

Next, using Lemma \ref{lem:useful-upb}, we have
\begin{align*}
|H_{5}^{[(n)]}| \leq \sum_{0 \leq k < \ell} 2^{-k} |\langle \mu, \Pp(\tilde{f}_{n}\, g_{k,\ell,n}) \rangle| \ind_{\{\ell \leq p\}} \leq C \,\sum_{\ell > 0} (2^{-\ell + 1} \, h_{n}^{2d} + \sum_{k = 0}^{\ell - 2} 2^{-k} \, \alpha^{\ell - k} \, h_{n}^{3d/2}) \leq C \, h_{n}^{3d/2}.
\end{align*}
This implies that $\lim_{n \rightarrow \infty} H_{5}^{[(n)]} = 0$  and then that $\lim_{n \rightarrow \infty} V_{5}(n) = 0$ a.s. and in $L^{2}.$

\medskip
 
Finally, since $V_{2}(n) = V_{5}(n) + V_{6}(n),$ it follows from the foregoing that $\lim_{n \rightarrow \infty} \EE[V_{2}(n)^{2}] = 0$ and this ends the proof.
\end{proof}
Now we treat the term $V_{1}(n)$. Recall $$V_{1}(n) = |\GG_{n}|^{-1}\sum_{i \in \GG_{n-p}}\sum_{k = 0}^{p} \EE_{X_{i}}[M_{\GG_{p - k}}(\tilde{f}_{n})^{2}].$$ We have the following convergence.

\begin{lem}\label{lem:cvV1est}
Under the assumptions of Theorem \ref{thm:flx-T}, we have the following convergence:
\begin{equation*}
\lim_{n\rightarrow \infty} V_{1}(n) = 2 \, \|K_{0}\|_{2}^{6} \, \mu^{\vt}(x,x_{0},x_{1}) \quad \text{in probability}.
\end{equation*}
\end{lem}

\begin{proof}
Let $k \in \{0,\ldots,p\}$ and $i \in \GG_{n-p}$. Conditioning by $\Hh_{i,p-k}$ and using the branching Markov property, we get
\begin{equation*}
\EE_{X_{i}}[M_{\GG_{p - k}}(\tilde{f}_{n})^{2}] \, = \, \EE_{X_{i}}[M_{\GG_{p-k}}(\Pp(\tilde{f}_{n}^{2}) - (\Pp \tilde{f}_{n})^{2})]  \, + \, \EE_{X_{i}}[(M_{\GG_{p-k}}(\Pp \tilde{f}_{n}))^{2}].
\end{equation*}
Using the latter inequality and the fact that $\Pp(\tilde{f}_{n}^{2}) - (\Pp \tilde{f}_{n})^{2} = \Pp(f_{n}^{2}) - (\Pp f_{n})^{2},$ we get
\begin{equation*}
V_{1}(n) = V_{3}(n) + V_{4}(n) - V_{7}(n),
\end{equation*}
where
\begin{align*}
V_{3}(n) \, = \, |\GG_{n}|^{-1} \, \sum_{i \in \GG_{n-p}} \, \sum_{k = 0}^{p} \EE_{X_{i}}[M_{\GG_{p-k}}(\Pp(f_{n}^{2}))]; \\ V_{7}(n) \, = \, |\GG_{n}|^{-1} \, \sum_{i \in \GG_{n-p}} \sum_{k = 0}^{p} \EE_{X_{i}}[M_{\GG_{p-k}}((\Pp f_{n})^{2})]; \\ V_{4}(n) \, = \, |\GG_{n}|^{-1} \, \sum_{i \in \GG_{n-p}} \sum_{k = 0}^{p} \EE_{X_{i}}[(M_{\GG_{p-k}}(\Pp \tilde{f}_{n}))^{2}].
\end{align*}
First we treat $V_{7}(n)$. We set 
\begin{equation*}
H_{7,n} = \sum_{k = 0}^{p} 2^{-k} \, \Qq^{p-k} ((\Pp f_{n})^{2}).
\end{equation*}
Using \eqref{eq:Q1}, we have $V_{7}(n) = |\GG_{n-p}|^{-1} M_{\GG_{n-p}}(H_{7,n}).$  Using Lemma \ref{lem:useful-upb} and distinguishing the cases $k = p$ and $k \leq p-1$, we get $|V_{7}(n)| \leq \, C \, (2^{-p} \,h_{n}^{d} \, + \, h_{n}^{2d}).$ It then follows that $\lim_{n \rightarrow \infty} V_{7}(n) = 0$ in probability.
 
\medskip

Next, we treat the term $V_{3}(n)$. We set $A_{3,n} = H_{3,n} - H_{3}^{[n]}$, with:
\begin{equation}\label{eq:H3crn-T}
H_{3,n} = \sum_{k = 0}^{p} 2^{-k} \, \Qq^{p-k}(\Pp( f_{n}^{2})) \quad \text{and} \quad H_{3}^{[n]} = \sum_{k=0}^{p} 2^{-k} \, \langle \mu,\Pp(f_{n}^{2}) \rangle = 2(1 - 2^{-p-1}) \langle \mu,\Pp(f_{n}^{2}) \rangle.
\end{equation} 
We set $g_{n} = \Pp(f_{n}^{2}) - \langle \mu,\Pp(f_{n}^{2}) \rangle.$   Using \eqref{eq:geom-ergB} and Lemma \ref{lem:useful-upb}, we have
\begin{align}
|V_{3}(n) - H_{3}^{[n]}| \, &\leq |\GG_{n-p}|^{-1} \, M_{\GG_{n-p}}(2^{-p} \, |g_{n}| )+ M_{\GG_{n-p}}(\sum_{k = 0}^{p-1} 2^{-k} |\Qq^{p-k-1}(\Qq(g_{n}))|) \nonumber \\
& \leq \, C \, 2^{-p} \, + \, C \, |\GG_{n-p}|^{-1} \, M_{\GG_{n-p}}(\sum_{k = 0}^{p-1} 2^{-k} \, \alpha^{p-k} \, \|\Qq(\Pp(f_{n}^{2}))\|_{\infty} \, V) \nonumber \\
& \leq \, C \, 2^{-p} \, + \, C \, a_{n} \, |\GG_{n-p}|^{-1} M_{\GG_{n-p}}(V), \label{eq:I-V3n-H3cn}
\end{align}
where
\begin{equation*}
a_{n} = \begin{cases} 2^{-p} & \text{if $2\alpha < 1$} \\ p\, 2^{-p} & \text{if $2\alpha = 1$} \\ \alpha^{p} & \text{if $2\alpha > 1.$}\end{cases}
\end{equation*}
Using \eqref{eq:lfgn-G} and the fact that $\lim_{n \rightarrow +\infty} a_{n} = 0$, we find that
\begin{equation*}
\lim_{n \rightarrow +\infty} C \, 2^{-p} \, + \, C \, a_{n} \, |\GG_{n-p}|^{-1} M_{\GG_{n-p}}(V) = 0 \quad \text{a.s. and in $L^{2}.$}
\end{equation*}
From \eqref{eq:I-V3n-H3cn}, this implies that
\begin{equation*}
\lim_{n \rightarrow + \infty}|V_{3}(n) - H_{3}^{[n]}| = 0 \quad \text{in probability.}
\end{equation*}  
 
Using Lemma \ref{lem:bochner}, we get
\begin{equation*}
\lim_{n \rightarrow \infty} H_{3}^{[n]} = \lim_{n \rightarrow \infty} 2 \, (1 - 2^{-p-1}) \, \langle \mu,\Pp(f_{n}^{2}) =  2 \, \|K_{0}\|_{2}^{6} \, \mu^{\vt}(x,x_{0},x_{1}).
\end{equation*}
From the foregoing, we conclude that $\lim_{n \rightarrow \infty} V_{3}(n) = 2 \, \|K_{0}\|_{2}^{6} \, \mu^{\vt}(x,x_{0},x_{1})$ in probability.

\medskip

Finally, we treat the term $V_{4}(n)$. Using \eqref{eq:Q2}, we have
\begin{equation*}
V_{4}(n) \, = \, V_{8}(n) \, + \, V_{9}(n),
\end{equation*}
where
\begin{equation*}
V_{8}(n) \, = \, |\GG_{n-p}|^{-1} \, M_{\GG_{n-p}}(H_{8,n}) \quad \text{and} \quad V_{9}(n) \, = \, |\GG_{n-p}|^{-1} \, M_{\GG_{n-p}}(H_{9,n}), 
\end{equation*}
with
\begin{equation*}
H_{8,n} = \sum_{k = 0}^{p} 2^{-k} \Qq^{p-k} ((\Pp \tilde{f}_{n})^{2}) \quad \text{and} \quad H_{9,n} = \sum_{\substack{k \geq 0, \ell \geq 0}} 2^{- k + \ell} \Qq^{p - k - \ell - 1} \Pp(\Qq^{\ell} \Pp \tilde{f}_{n} \otimes \Qq^{\ell} \Pp \tilde{f}_{n}) \ind_{\{k + \ell < p\}}.
\end{equation*}
Writing
\begin{equation*}
H_{8}^{[n]} = \sum_{k = 0}^{p} 2^{-k} \langle \mu, (\Pp \tilde{f}_{n})^{2}\rangle \quad \text{and} \quad H_{9}^{[n]} = \sum_{\substack{k \geq 0, \ell \geq 0}} 2^{- k + \ell} \langle \mu,\Pp(\Qq^{\ell} \Pp \tilde{f}_{n} \otimes \Qq^{\ell} \Pp \tilde{f}_{n})\rangle \ind_{\{k + \ell < p\}},
\end{equation*}
we prove, as previously, that
\begin{align*}
\lim_{n \rightarrow \infty} |V_{8}(n) - H_{8}^{[n]}| = \lim_{n \rightarrow \infty} |V_{9}(n) - H_{9}^{[n]}| = 0 \quad \text{a.s. and in $L^{2}$;} \quad
\lim_{n \rightarrow \infty} H_{8}^{[n]} = \lim_{n \rightarrow \infty} H_{9}^{[n]}  = 0.
\end{align*}
As a result, we find that
\begin{equation*}
\lim_{n \rightarrow \infty} |V_{8}(n)| = \lim_{n \rightarrow \infty} |V_{9}(n)| = 0 \quad \text{in probability.}
\end{equation*}
Since $V_{4}(n) = V_{8}(n) + V_{9}(n)$, we conclude that $\lim_{n \rightarrow \infty} V_{4}(n) = 0$ in probability. Finally, since $V_{1}(n) = V_{3}(n) + V_{4}(n) - V_{7}(n)$, the result of the Lemma follows from the foregoing.
\end{proof}

As a consequence of \eqref{eq:dVtriang}, Lemmas \ref{lem:cvR2est}, \ref{lem:cvV2est} and \ref{lem:cvV1est}, we have the following result.

\begin{lem}\label{lem:cvV-est-T}
Under the assumptions of Theorem \ref{thm:flx-T}, we have the following convergence:
\begin{equation*}
\lim_{n\rightarrow \infty} V(n) = 2 \, \|K_{0}\|_{2}^{6} \, \mu^{\vt}(x,x_{0},x_{1}) \quad \text{in probability}.
\end{equation*}
\end{lem}

We now check the Lindeberg condition using a fourth moment condition. We set:
\begin{equation*}\label{eq:def-R3T}
R_3(n) = \sum_{i\in \G_{n-p_n}} \E\left[\Delta_{n,i}(f_{n})^4\right].
\end{equation*}

\begin{lem}\label{lem:cvR3T}
Under the assumptions of Theorem \ref{thm:flx-T}, we have that $\lim_{n \rightarrow \infty} R_{3}(n) = 0$
\end{lem}

\begin{proof}
We have
\begin{align}
R_{3}(n) &\leq 16 \, (p + 1)^{3} |\GG_{n}|^{-2} \sum_{i \in \GG_{n-p}}\sum_{\ell = 0}^{p} \EE[(M_{i\GG_{p-\ell}}(\tilde{f}_{n}))^{4}] \nonumber \\ & \leq 128 \, (p + 1)^{3} \, |\GG_{n}|^{-2} \sum_{i \in \GG_{n-p}} \sum_{\ell = 0}^{p} \EE[(M_{i\GG_{p-\ell}}(f_{n} - \Pp f_{n}))^{4}] \label{eq:BR3T} \\ & \hspace{2cm} + \, 128 \, (p +  1)^{3} \, |\GG_{n}|^{-2} \, \sum_{i \in \GG_{n-p}} \sum_{\ell = 0}^{p} \EE[(M_{i\GG_{p-\ell}}(\Pp \tilde{f}_{n}))^{4}], \nonumber
\end{align}
where we used that $(\sum_{k=0}^r a_k)^4 \leq  (r+1)^3 \sum_{k=0}^r a_k^4$ for the two inequalities (resp. with $r=1$ and $r=p$), Jensen inequality for the first inequality and the decomposition $f_{n} = (f_{n} - \Pp f_{n}) + \Pp f_{n}$ for the last inequality. For the second term of the right hand side of \eqref{eq:BR3T}, we follow the proof of Lemma 5.6 in \cite{BD3} and Lemma 4.7 in \cite{BD1} to find that
\begin{equation*}
n^{3} \, |\GG_{n}|^{-2} \, \sum_{i \in \GG_{n-p}} \sum_{\ell = 0}^{p} \EE[(M_{i\GG_{p-\ell}}(\Pp \tilde{f}_{n}))^{4}] \leq C \, n^{5} (2^{-n+p} h_{n} \ind_{\{2\alpha^{2} \leq 1\}} + 2^{-n + p} \, (2\alpha^{2})^{2p} h_{n}^{6d} \ind_{\{2\alpha^{2} > 1\}}),
\end{equation*} 
and using \eqref{eq:2alpha2h3d}, this implies that
\begin{equation}\label{eq:cvR32ntri}
\lim_{n \rightarrow +\infty} n^{3} \, |\GG_{n}|^{-2} \, \sum_{i \in \GG_{n-p}} \sum_{\ell = 0}^{p} \EE[(M_{i\GG_{p-\ell}}(\Pp \tilde{f}_{n}))^{4}] = 0.
\end{equation}  
We are now going to treat the first term of \eqref{eq:BR3T}. Since $\Pp(f_{n} - \Pp(f_{n})) = 0$, we have, (see Remark 2.3 in \cite{BDG14} for more details),
\begin{equation}\label{eq:R23BDG}
\EE_{x}[(M_{\GG_{p-\ell}}(f_{n} - \Pp f_{n}))^{4}] \leq \, g_{n,\ell}(x) \, + 6 \, h_{n,\ell}(x),
\end{equation}
with:
\begin{equation*}
g_{n,\ell}(x) = \EE_{x}[M_{\GG_{p-\ell}}(\Pp((f_{n} - \Pp f_{n})^{4}))] \quad \text{and} \quad h_{n,\ell}(x) = \EE_{x}[(M_{\GG_{p-\ell}}(\Pp((f_{n} - \Pp f_{n})^{2})))^{2}].
\end{equation*}
We set
\begin{equation*}
R_{3,1}(n) = (p + 1)^{3} \, 2^{-2n} \sum_{i \in \GG_{n-p}} \sum_{\ell = 0}^{p} \EE[(M_{i\GG_{p-\ell}}(f_{n} - \Pp f_{n}))^{4}].
\end{equation*}
Using the branching Markov property and \eqref{eq:R23BDG} for the first inequality and \eqref{eq:Q1} for equality, we get
\begin{align}
\nonumber R_{3,1}(n) &\leq C \, n^{3} \, 2^{- 2n} \, \sum_{\ell = 0}^{p} \EE[M_{\GG_{n - p}}(g_{n,\ell})] + C \, n^{3} \, 2^{- 2n} \, \sum_{\ell = 0}^{p} \EE[M_{\GG_{n-p}}(h_{n,\ell})]
\\ 
\label{eq:IR31ntri} & = \, C \, n^{3} \, 2^{-n-p} \, \sum_{\ell = 0}^{p} \langle \nu, \Qq^{n-p} g_{n,\ell} \rangle \, + \, C \, n^{3} \, 2^{-n-p} \, \sum_{\ell = 0}^{p} \langle \nu, \Qq^{n-p} h_{n,\ell} \rangle.
\end{align}
Using Lemma \ref{lem:useful-upb}, we get
\begin{equation*}
\langle \nu, \Qq^{n-p} g_{n,\ell} \rangle \, = \, 2 ^{p-\ell} \, \langle \nu, \Qq^{n-\ell}(\Pp((f_{n} - \Pp f_{n})^{4})) \rangle \, \leq  \, 2^{p-\ell} \, \|\Qq\Pp((f_{n} - \Pp f_{n})^{4})\|_{\infty} \leq \, C \, h_{n}^{-3d} \, 2^{p - \ell}.
\end{equation*}
The latter inequality implies that
\begin{equation}\label{eq:cvR31nP1}
n^{3} \, 2^{-n-p} \, \sum_{\ell = 0}^{p} \langle \nu, \Qq^{n-p} g_{n,\ell} \rangle  \leq C \, n^{3} \, (2^{n}h_{n}^{3d})^{-1}.
\end{equation} 
Using \eqref{eq:Q2}, the fact that $\Pp((f_{n} - \Pp f_{n})^{2}) \leq \Pp (f_{n}^{2})$ for the first inequality and Lemma \ref{lem:useful-upb} for the second inequality, we get
\begin{align*}
\langle \nu, \Qq^{n-p} h_{n,\ell}\rangle & \leq 2^{p-\ell} \, \langle \nu,\Qq^{n-\ell}(\Pp (f_{n}^{2}))^{2} \rangle + \sum_{k = 0}^{p - \ell - 1} 2^{p-\ell+k} \langle \nu, \Qq^{n-\ell-k-1}\Pp(\Qq^{k}\Pp(f_{n}^{2} \otimes^{2}) \rangle
\\ & \leq \, C \,  2^{p-\ell} \, + \, C \, 2^{2(p-\ell)}.
\end{align*}
The latter inequality implies that
\begin{equation}\label{eq:cvR31nP2}
n^{3} \, 2^{-n-p} \, \sum_{\ell = 0}^{p} \langle \nu, \Qq^{n-p} h_{n,\ell} \rangle \, \leq \, C \, n^{3}(2^{-n} \, + \, 2^{- n + p}).
\end{equation}
From \eqref{eq:IR31ntri}, \eqref{eq:cvR31nP1} and \eqref{eq:cvR31nP2}, we conclude that $\lim_{n \rightarrow \infty} R_{3,1}(n) = 0.$ Finally, from \eqref{eq:BR3T} and \eqref{eq:cvR32ntri}, this proves that $\lim_{n \rightarrow \infty} R_{3}(n) = 0$.
\end{proof}

We can  now use  Theorem 3.2 and  Corollary 3.1, p.~58, and  the remark p.~59  from  \cite{hh:ml}  to  deduce  from  Lemmas \ref{lem:cvV-est-T}  and \ref{lem:cvR3T} that  $\Delta_n(f_{n})$ converges in distribution  towards a Gaussian  real-valued   random  variable  with   deterministic  variance $\sigma^{2}$.  Using \eqref{eq:DNn0k0-T} and Lemmas  \ref{lem:cge-R0-est}  and  \ref{lem:cge-R1-est}, we  then  deduce  Theorem \ref{thm:flx-T} for $A_{n} = \TT_{n}$.

\medskip

For $A_{n} = \GG_{n}$, we have
\begin{equation*}
N_{n,\emptyset}(f_{n}) = R_{1}(n) + \Delta_{n}(f_{n}),
\end{equation*}
where
\begin{equation*}
R_{1}(n) = \sum_{i \in \GG_{n - p}} \EE[N_{n,i}(f_{n})|\Ff_{i}] \quad \text{and} \quad \Delta_{n}(f_{n}) =  \sum_{i \in \GG_{n-p}} \Delta_{n,i}(f_{n}),
\end{equation*}
and for all $i \in \GG_{n-p}$,
\begin{equation*}
N_{n,i}(f_{n}) = M_{i\GG_{p}}(\tilde{f}_{n}) \quad \text{and} \quad \Delta_{n,i}(f_{n}) = N_{n,i}(f_{n}) - \EE[N_{n,i}(f_{n})|\Ff_{i}].
\end{equation*}
Following exactly the proof of Lemma \ref{lem:cge-R1-est}, \ref{lem:cvV-est-T}  and \ref{lem:cvR3T}, we get the result for this case. We note that for $\A_{n} = \GG_{n}$, the factor $2$ is missing in the asymptotic variance. This come from the fact that here, $H_{3}^{[n]}$ defined in \eqref{eq:H3crn-T} is simply equal to $\langle \mu,\Pp(f_{n}^{2}) \rangle$. 

\section{Proof of Theorem \ref{thm:cltEst-T}}\label{proof:thm:cltEst-T}

We begin the proof with $\A_{n} = \TT_{n}$. From \eqref{eq:Dmuhatvt-T}, we have
\begin{align*}
|\TT_{n}|^{1/2} \, h^{3d/2} \, (\mu_{\TT_{n}}^{\vt}(xx_{0}x_{1}) - \mu^{\vt}(xx_{0}x_{1})) & = (|\GG_{n}|/|\TT_{n}|)^{1/2} \, N_{n,\emptyset}(f_{n})  \, + \, B_{n}(xx_{0}x_{1}),
\end{align*}
where the bias term $B_{n}(xx_{0}x_{1})$ is defined by
\begin{equation*}
B_{n}(xx_{0}x_{1}) = |\TT_{n}|^{1/2} \, h_{n}^{3d/2} \, (h^{-3/2} \, \langle \mu^{\vt}, f_{n} \rangle - \mu^{\vt}(xx_{0}x_{1})).
\end{equation*}
Since $\lim_{n \rightarrow \infty} (|\GG_{n}|/|\TT_{n}|)^{1/2} = 1/\sqrt{2}$, from Theorem \ref{thm:flx-T}, it suffices, to obtain the result of Theorem \ref{thm:cltEst-T}, to prove that $\lim_{n \rightarrow \infty} B_{n}(xx_{0}x_{1}) = 0.$ 
Using the Taylor expansion and Assumption \ref{hyp:estim-tcl}, one can prove that (see \cite{BD2020} for more details) 
\begin{equation*}
B_{n}(xx_{0}x_{1}) \, \leq \, C \, \sqrt{|\TT_{n}|h_{n}^{2s + 3d}}.
\end{equation*}
Since $\lim_{n \rightarrow \infty} |\TT_{n}|h_{n}^{2s + 3d} = 0$, we conclude that $\lim_{n \rightarrow \infty} B_{n}(xx_{0}x_{1}) = 0$ and this ends the proof for $\A_{n} = \TT_{n}.$

\medskip

For $\A_{n} = \GG_{n}$ the proof follows exactly the same lines.

\section{Proof of Theorem \ref{thm:cltPhat-sub}}\label{sec:cltphat-sub}

First of all, we have the following decomposition:
\begin{align*}
\sqrt{|\A_{n}|h_{n}^{3d}} (\widehat{\Pp}_{\A_{n}}(xx_{0}x_{1}) - \Pp(xx_{0}x_{1})) & = (|\A_{n}|h_{n}^{3d})^{1/2} \, (\frac{\widehat{\mu}^{\vt}_{\A_{n}}(xx_{0}x_{1})}{\widehat{\mu}_{\A_{n}}(x)} - \frac{\mu^{\vt}(xx_{0}x_{1})}{\widehat{\mu}_{\A_{n}}(x)}) 
\\  & \hspace{2cm} - \, \frac{\mu^{\vt}(xx_{0}x_{1})}{\mu(x)\widehat{\mu}_{\A_{n}}(x)} \, (|\A_{n}|h_{n}^{3d})^{1/2} \,(\widehat{\mu}_{\A_{n}}(x) - \mu(x)).
\end{align*}
Then, the proof of Theorem \ref{thm:cltPhat-sub} is a direct consequence of the previous decomposition and Lemmas \ref{eq:LprPhat-P-R} and \ref{eq:LmainPhat-P} below.

\begin{lem}\label{eq:LprPhat-P-R}
Under Assumptions of Theorem \ref{thm:cltPhat-sub}, we have
\begin{equation*}
\lim_{n \rightarrow \infty} \frac{\mu^{\vt}(xx_{0}x_{1})}{\mu(x)\widehat{\mu}_{\A_{n}}(x)} \, (|\A_{n}|h_{n}^{3d})^{1/2} \,(\widehat{\mu}_{\A_{n}}(x) - \mu(x)) = 0 \quad \text{in probability.}
\end{equation*}
\end{lem}

\begin{proof}
We consider the function $g_{n}$ defined on $S$ by: $g_{n}(y) = h_{n}^{-d/2} K_{0}(h_{n}^{-1}(x-y))$ for all $y \in S$. We begin the proof with $\A_{n} = \TT_{n}.$ We set $\tilde{g}_{n} = g_{n} - \langle \mu, g_{n} \rangle.$ We have the following decomposition:
\begin{equation}\label{eq:D-muhat-mu}
\widehat{\mu}_{\TT_{n}}(x) - \mu(x) = |\TT_{n}|^{-1} h_{n}^{-d/2} (\sum_{\ell = 0}^{2} M_{\GG_{\ell}}(\tilde{g}_{n}) + \sum_{\ell = 3}^{n} M_{\GG_{\ell}}(\tilde{g}_{n})) + \langle \mu,h_{n}^{-d/2} g_{n} \rangle - \mu(x).  
\end{equation}
Using the fact that $K_{0}$ is bounded, integration by parts and Assumption \ref{hyp:ub-density}, we have the following upper bounds: 
\begin{equation}\label{eq:ub-gn-Qmu}
\|g_{n}\|_{\infty} \leq \|K_{0}\| _{\infty} \, h_{n}^{-d/2}; \quad \|\Qq g_{n}\|_{\infty} + |\langle \mu,g_{n} \rangle| \leq 2 \, C_{0} \, \|K_{0}\|_{1} \, h_{n}^{d/2}.
\end{equation}
Using \eqref{eq:ub-gn-Qmu}, we find that
\begin{equation}\label{eq:Lubmuhatmu1}
|\TT_{n}|^{-1} h_{n}^{-d/2} |\sum_{\ell = 0}^{2} M_{\GG_{\ell}}(\tilde{g}_{n})| \leq C \, (|\TT_{n}| \, h_{n}^{d})^{-1}.
\end{equation}
Next, from Minkowski's inequality, we have
\begin{equation*}
\EE[(\sum_{\ell = 3}^{n} M_{\GG_{\ell}})^{2}] \leq (\sum_{\ell = 3}^{n} (\EE[(M_{\GG_{\ell}}(\tilde{g}_{n}))^{2}])^{1/2})^{2}.
\end{equation*}
Using  \eqref{eq:Q2}, \eqref{eq:geom-ergB} and \eqref{eq:ub-gn-Qmu}, we get:
\begin{align*}
\EE[(M_{\GG_{\ell}}(\tilde{g}_{n}))^{2}] &\leq C \, 2^{\ell} \, \langle \nu, \Qq^{\ell}(\tilde{g}_{n}^{2}) \rangle \, + \, \sum_{r = 0}^{\ell - 1} 2^{\ell + r} \langle \nu, \Qq^{\ell - r - 1}(\Pp(|\Qq^{r} \tilde{g}_{n}| \otimes^{2})) \rangle \\
& \leq \, C \,2^{\ell} \, + \, C \, h_{n}^{d} \, \sum_{r=0}^{\ell - 1} 2^{\ell} \, (2\alpha^{2})^{r} 
\leq \, C \, 2^{\ell} \, \ind_{\{2\alpha^{2} \leq 1\}}  \, + \,  C \, 2^{\ell}  (1  +  (2\alpha^{2})^{\ell} \, h_{n}^{d}) \ind_{\{2\alpha^{2} > 1\}}.
\end{align*}
The latter inequality implies that
\begin{equation}\label{eq:Lubmuhatmu2}
\EE[(|\TT_{n}|^{-1} h_{n}^{-d/2} \sum_{\ell = 3}^{n} M_{\GG_{\ell}}(\tilde{g}_{n}))^{2}] \leq C \, (|\TT_{n}| h_{n}^{3d})^{-1} \, (h_{n}^{2d} + (2\alpha^{2})^{n} h_{n}^{3d} \, \ind_{\{2\alpha^{2} > 1\}}).
\end{equation}
Using \eqref{eq:Lubmuhatmu1}, \eqref{eq:Lubmuhatmu2} and \eqref{eq:2alpha2h3d}, we deduce that
\begin{equation}\label{eq:lim-E-Tn-hd}
\lim_{n \rightarrow +\infty} \EE[(|\TT_{n}|^{-1} h_{n}^{-d/2} \sum_{\ell = 0}^{n} M_{\GG_{\ell}}(\tilde{g}_{n}))^{2}] = 0.
\end{equation}
 
Next, using Taylor expansion and Assumption \ref{hyp:estim-tcl}, we get (see \cite{BD2020} for more details)
\begin{equation}\label{eq:Lubmuhatmu3}
|\langle \mu,h_{n}^{-d/2} g_{n} \rangle - \mu(x)| \leq C \, h_{n}^{s}.
\end{equation}
From \eqref{eq:D-muhat-mu}, \eqref{eq:lim-E-Tn-hd} and \eqref{eq:Lubmuhatmu3}, we deduce that
\begin{equation}\label{eq:LSTnmuhat-mu} 
\lim_{n \rightarrow \infty} \widehat{\mu}_{\TT_{n}}(x) = \mu(x) \quad \text{in probability.} 
\end{equation}
We further deduce that
\begin{equation*}
\EE[|\TT_{n}|h_{n}^{3d} (\widehat{\mu}_{\TT_{n}}(x) - \mu(x))^{2}] \, \leq \, C \, (h_{n}^{2d} \, + \, (2\alpha^{2})^{n} h_{n}^{3d} \ind_{\{2\alpha^{2} > 1\}} \, + \, |\TT_{n}| h_{n}^{3d + 2s}).
\end{equation*}
Using \eqref{eq:2alpha2h3d}, the latter inequality implies that
\begin{equation}\label{eq:LTTnmuhat-mu}
\lim_{n \rightarrow \infty} (|\TT_{n}|h_{n}^{3d})^{1/2} \,(\widehat{\mu}_{\TT_{n}}(x) - \mu(x)) = 0 \quad \text{in probability.}
\end{equation}
From \eqref{eq:LSTnmuhat-mu}, \eqref{eq:LTTnmuhat-mu} and using Slutsky's Lemma, we get 
\begin{equation*}
\lim_{n \rightarrow \infty} \frac{\mu^{\vt}(xx_{0}x_{1})}{\mu(x)\widehat{\mu}_{\TT_{n}}(x)} \, (|\TT_{n}|h_{n}^{3d})^{1/2} \,(\widehat{\mu}_{\TT_{n}}(x) - \mu(x)) = 0 \quad \text{in probability.}
\end{equation*}
For $\A_{n} = \GG_{n}$, we follows exactly the same lines and this ends the proof.
\end{proof}

\begin{lem}\label{eq:LmainPhat-P}
Under Assumptions of Theorem \ref{thm:cltPhat-sub}, we have
\begin{equation*}
(|\A_{n}|h_{n}^{3d})^{1/2} \, (\frac{\widehat{\mu}^{\vt}_{\A_{n}}(xx_{0}x_{1})}{\widehat{\mu}_{\A_{n}}(x)} - \frac{\mu^{\vt}(xx_{0}x_{1})}{\widehat{\mu}_{\A_{n}}(x)}) \inL G,
\end{equation*}
where $G$ is a centered Gaussian real-valued random variable with mean $0$ and variance $$\sigma^{2} \, = \, \|K_{0}\|_{2}^{6}\, \Pp(x,x_{0},x_{1})/\mu(x).$$ 
\end{lem}

\begin{proof}
This is a direct consequence of Theorem \ref{thm:cltEst-T}, \eqref{eq:LSTnmuhat-mu} and Slutsky's Lemma. 
\end{proof}

\section{Proof of \eqref{eq:AMSE}}\label{sec:proof-amise}

We set $f_{0h}(y) = h^{-1} K_{0}(h^{-1}(x-y))$ and recall $\langle \mu, f_{0h} \rangle = \int_{\RR}f_{0h}(y)\mu(y)dy.$ Using the decomposition
\begin{equation*}
\widehat{\mu}(x) - \mu(x) = \frac{1}{|\GG_{n}|} \sum_{u \in \GG_{n}} \tilde{f}_{0h}(X_{u}) \, +  \langle \mu, f_{0h} \rangle \, - \, \mu(x),
\end{equation*}
we obtain the following biais-variance type decomposition.
\begin{equation}\label{eq:Decom-BAR}
\EE\left[(\widehat{\mu}(x) - \mu(x))^{2}\right] \leq 2 (|\GG_{n}|)^{2} \EE\big[( \sum_{u \in \GG_{n}} \tilde{f}_{0h}(X_{u}))^{2}\big] + 2 \big(\langle \mu, f_{0h} \rangle \, - \, \mu(x)\big)^{2}.
\end{equation}
Using \eqref{eq:Q2} and the fact that the process $\Ll(X_{\emptyset}) = \mu$ (which implies that $\mu \Qq = \mu$), we get
\begin{equation}\label{eq:var-BAR}
|\GG_{n}|^{-2} \, \EE\big[( \sum_{u \in \GG_{n}} \tilde{f}_{0h}(X_{u}))^{2}\big] \, \leq \, \frac{2 \, \langle \mu, \tilde{f}_{0h}^{2} \rangle}{|\GG_{n}|} \, + \, \frac{1}{|\GG_{n}|} \sum_{k=1}^{n-1} 2^{k} \langle \mu, (\Qq^{k-1}(\Qq \tilde{f}_{0h}))^{2} \rangle.
\end{equation}
We now plan to use \eqref{eq:erg-BAR} with $f = \Qq f_{0h}$. For all $y \in \RR,$ we get, after the change of variable $t = h^{-1}(x-z)$ and the use of the first-order Taylor's expansion,
\begin{align*}
(\Qq f_{0h})'(y) &= \frac{1}{h} \int_{\RR} K_{0}(h^{-1}(x-z)) \frac{\partial \Qq}{\partial y}(y,z) dz \\ 
& = \int_{\RR} K_{0}(z) \frac{\partial \Qq}{\partial y}(y, x \, - \, hz) \, dz \\
& = \frac{\partial \Qq}{\partial y}(y,x) \int_{\RR} K_{0}(z) \, dz + \smallO(1) = \frac{\partial \Qq}{\partial y}(y,x) \, + \, \smallO(1).
\end{align*}
We then have that
\begin{equation*}
\|(\Qq f_{0h})'\|_{\infty} = \sup_{y \in \RR} \Big\{\frac{\partial \Qq}{\partial y}(\cdot,x)\Big\} \, + \, \smallO(1) \, = \, \frac{1}{\sqrt{2\pi} \sigma^{2}} \expp{-1/2} + \smallO(1).
\end{equation*}
Using the latter equality and \eqref{eq:erg-BAR}, we get, for all $k \geq 1$,
\begin{equation}\label{eq:A-lyapunov}
|\Qq^{k-1}(\Qq \tilde{f}_{0h})|(y) \leq  \frac{1}{\sqrt{2\pi} \sigma^{2} \expp{1/2} (1-a)} (\sigma(1 + a)^{-1} + |y|) a^{k-1} + \smallO(1).
\end{equation}
Recall $\mu$ is the Gaussian law $\Nn(0,\sigma^{2}(1-a^{2}))$. Using \eqref{eq:A-lyapunov} and $(\sigma(1 + a)^{-1} + |y|)^{2} \leq 2 \, \sigma^{2}(1 + a)^{-2} \, + \, 2 \, y^{2}$, we get
\begin{equation}\label{eq:A-lya-C}
\langle \mu, (\Qq^{k-1}(\Qq \tilde{f}_{0h}))^{2} \rangle \leq  \frac{1}{e \pi \sigma^{2}} \Big(\frac{1}{(1 - a)^{2}} + \frac{1}{1 - a^{2}}\Big) a^{2(k-1)} \, + \, \smallO(1) \, \leq \, \frac{C_{a,\sigma}}{a^{2}} \, a^{2k} \, + \, \smallO(1).
\end{equation}
Now, \eqref{eq:A-lya-C} and \eqref{eq:var-BAR} implies that
\begin{equation*}
(|\GG_{n}|)^{2} \EE\big[( \sum_{u \in \GG_{n}} \tilde{f}_{0h}(X_{u}))^{2}\big] \, \leq \, \frac{2 \, \langle \mu, \tilde{f}_{0h}^{2} \rangle}{|\GG_{n}|} \, + \, \frac{M_{a,\sigma}}{a^{2} |\GG_{n}|} \sum_{k=1}^{n-1} (2a^{2})^{k} \, + \, \smallO(1). 
\end{equation*}
Putting the latter inequality into \eqref{eq:Decom-BAR}, we obtain
\begin{equation}\label{eq:Decom-BAR+}
\EE\left[(\widehat{\mu}(x) - \mu(x))^{2}\right] \leq \, \frac{2 M_{a,\sigma}}{a^{2} |\GG_{n}|} \sum_{k=1}^{n-1} (2a^{2})^{k} \, + \, \frac{4 \, \langle \mu, \tilde{f}_{0h}^{2} \rangle}{|\GG_{n}|} \, +   \, 2 \, \big(\langle \mu, f_{0h} \rangle \, - \, \mu(x)\big)^{2} \, + \, \smallO(1).
\end{equation}
Finally, it is very standard to get asymptotic equivalence of the second and the third term of the right hand side of \eqref{eq:Decom-BAR+} (see for e.g. \cite{silverman1986density}, Section 3.3.1 for more details). This ends the proof of \eqref{eq:AMSE}.

\section{Appendix}\label{sec:appendix}

First, we give some useful upper bounds. We recall that $S = \RR^{d}$. Recall $f_{n}$ defined in \eqref{eq:def-fnT}.
\begin{lem}\label{lem:useful-upb}
Under Assumption \eqref{hyp:ub-density}, we have:
\begin{align*}\label{eq:EPpfn2}
&\|\Pp |f_{n}|\|_{\infty} \leq \|K_{0}\|_{1}^{2} \, \|K_{0}\|_{\infty} \, \|\Pp\|_{\infty} \, h_{n}^{d/2};
\\ &\|\Qq \Pp |f_{n}|\|_{\infty} \leq \|\Pp\|_{\infty} \, \|\Qq\|_{\infty} \, \|K_{0}\|_{1}^{3} \, h_{n}^{3d/2};
\\ &|\langle \mu,\Pp f_{n}  \rangle| \leq \|\Pp\|_{\infty} \, \|\mu\|_{\infty} \, \|K_{0}\|_{1}^{3} \, h_{n}^{3d/2}
\\ &\langle \mu, \Pp (f_{n}^{2}) \rangle \leq \|\mu\|_{\infty} \, \|\Pp\|_{\infty} \, \|K_{0}\|_{2}^{6} 
\\ &\|\Qq \Pp f_{n}^{2}\|_{\infty} \leq \|K_{0}\|_{2}^{6} \, \|\Pp\|_{\infty} \, \|\Qq\|_{\infty};
\\ &\|\Pp(\Pp f_{n} \otimes \Pp f_{n})\|_{\infty} \leq \|\Pp\|_{\infty}^{3} \, \|K_{0}\|^{6}_{1} \, h_{n}^{3d};
\\ &\|\Qq((\Pp f_{n})^{2})\|_{\infty} \leq \|K_{0}\|_{2}^{2} \, \|K_{0}\|_{1}^{4} \, \|\Qq\|_{\infty} \, \|\Pp\|^{2}_{\infty} \, h^{2d};
\\ &\|\Pp((\Pp f_{n})^{2} \otimes (\Pp f_{n})^{2})\|_{\infty} \leq \|K_{0}\|_{2}^{4} \, \|K_{0}\|_{1}^{8} \, \|\Pp\|^{5}_{\infty} \, h^{4d};
\\ & \|\Qq(\Pp(f_{n}^{4}))\|_{\infty} \leq \|\Pp\|_{\infty} \, \|\Qq\|_{\infty} \, \|K_{0}\|_{4}^{12} \, h_{n}^{-3d}.
\end{align*} 
\end{lem}

\begin{proof}
Using a change of variables, we have, for all $y \in S$:
\begin{equation}\label{eq:EPpfn}
\Pp |f_{n}|(y) = h^{d/2} \, |K_{0}(h^{-1}(x-y))| \, \int_{S^{2}} |K_{0}|(y_{0})\, |K_{0}|(y_{1}) \Pp(y,x_{0} - h\,y_{0},x_{1} - h \, y_{1}) dy_{0} \, dy_{1}.
\end{equation}
This implies that
\begin{equation*}\label{eq:EPpfn2}
\|\Pp f_{n}\|_{\infty} \leq \|K_{0}\|_{1}^{2} \, \|K\|_{\infty} \, \|\Pp\|_{\infty} \, h_{n}^{d/2}.
\end{equation*}
From \eqref{eq:EPpfn} and using again a change of variables, we have, for all $t \in S$:
\begin{equation*}\label{eq:EQqPpfn}
\Qq \Pp |f_{n}|(t) = h_{n}^{3/2} \, \int_{S^{3}} |K_{0}|(y)\, |K_{0}|(y_{0})\, |K_{0}|(y_{1})\, \Pp(x - h\, y, x_{0} - h\, y_{0}, x_{1} - h\, y_{1}) \, \Qq(t, x - h\, y) dy\, dy_{0}\, dy_{1}.
\end{equation*}
This implies that
\begin{equation}\label{eq:EQqPpfn2}
\|\Qq \Pp |f_{n}|\|_{\infty} \leq \|\Pp\|_{\infty} \, \|\Qq\|_{\infty} \, \|K_{0}\|_{1}^{3} \, h_{n}^{3d/2}.
\end{equation}
As for \eqref{eq:EQqPpfn2}, we have
\begin{equation*}
|\langle \mu,\Pp f_{n}  \rangle| \leq \|\Pp\|_{\infty} \, \|\mu\|_{\infty} \, \|K_{0}\|_{1}^{3} \, h_{n}^{3d/2}.
\end{equation*}
Now, following the same ideas, we easily get the others upper bounds.
\end{proof}

We recall  the following result due to Bochner
(see \cite[Theorem 1A]{Parzen1962}  which can be easily
extended to any dimension $d\geq 1$). 
\begin{lem}
  \label{lem:bochner}
Let $(h_{n},n\in\NN)$ be a sequence of positive numbers converging to $0$ as $n$ goes to infinity. Let $g: \RR^{d} \rightarrow \RR$ be a measurable function such that $\int_{\RR^{d}} |g(x)|dx < +\infty$. Let $f: \RR^{d} \rightarrow \RR$  be a measurable function such that  $\norm{f}_{\infty}<+\infty $, $\int_{\R^d} |f(y)|\, dy < + \infty$  and $\lim_{|x|\rightarrow +\infty} |x|f(x)=0$. Define
\begin{equation*}
g_{n}(x) = h_{n}^{-d}\int_{\RR^{d}} f(h_{n}^{-1}(x-y))g(y)dy.
\end{equation*}
Then, we have at every point $x$ of continuity of $g$,
\begin{equation*}
\lim_{n\rightarrow +\infty} g_{n}(x) = g(x)\int_{\RR} f(y)dy.
\end{equation*} 
\end{lem}

In this section, we recall useful results on BMC from Bitseki-Delmas \cite{BD2020}.

\begin{lem}\label{lem:Qi}
Let $f,g\in \cb(S)$, $x\in S$ and $n\geq m\geq 0$. Assuming that all the quantities below are well defined,  we have:
\begin{align} 
\label{eq:Q1} &\E_x\left[M_{\G_n}(f)\right] = |\G_n|\, \cq^n f(x)= 2^n\, \cq^n f(x) ,\\ 
\label{eq:Q2} &\E_x\left[M_{\G_n}(f)^2\right] = 2^n\, \cq^n (f^2) (x) + \sum_{k=0}^{n-1} 2^{n+k}\,   \cq^{n-k-1}\left( \cp\left(\cq^{k}f\otimes \cq^k f \right)\right) (x), \\
\label{eq:Q2-bis} &\E_x\left[M_{\G_n}(f)M_{\G_m}(g)\right] = 2^{n} \cq^{m} \left(g \cq^{n-m} f\right)(x)\\ \nonumber &\hspace{4cm} + \sum_{k=0}^{m-1} 2^{n+k}\, \cq^{m-k-1} \left(\cp\left(\cq^k g \sot \cq^{n-m+k} f\right) \right)(x). 
\end{align}
\end{lem}

\bibliographystyle{abbrv}
\bibliography{biblio}

\begin{thebibliography}{10}

\bibitem{BD2020}
S.~V. Bitseki~Penda and J.-F. Delmas.
\newblock Central limit theorem for bifurcating markov chains, 2020.

\bibitem{BD2}
S.~V. Bitseki~Penda and J.-F. Delmas.
\newblock {C}entral limit theorem for bifurcating {M}arkov chains under
  ${L}^{2}$ ergodic conditions.
\newblock {\em Advances in Applied Probability}, pages 1--33, 2022.

\bibitem{BD1}
S.~V. Bitseki~Penda and J.-F. Delmas.
\newblock Central limit theorem for bifurcating markov chains under pointwise
  ergodic conditions.
\newblock {\em The Annals of Applied Probability}, 32(5):3817--3849, 2022.

\bibitem{BD3}
S.~V. Bitseki~Penda and J.-F. Delmas.
\newblock Central limit theorem for kernel estimator of invariant density in
  bifurcating markov chains models.
\newblock {\em Journal of Theoretical Probability}, pages 1--35, 2022.

\bibitem{BDG14}
S.~V. Bitseki~Penda, H.~Djellout, and A.~Guillin.
\newblock Deviation inequalities, moderate deviations and some limit theorems
  for bifurcating {M}arkov chains with application.
\newblock {\em Ann. Appl. Probab.}, 24(1):235--291, 2014.

\bibitem{comte-marie2021}
F.~Comte and N.~Marie.
\newblock On a nadaraya-watson estimator with two bandwidths.
\newblock {\em Electronic Journal of Statistics}, 15(1):2566--2607, 2021.

\bibitem{duflo2013random}
M.~Duflo.
\newblock {\em Random iterative models}, volume~34.
\newblock Springer Science \& Business Media, 2013.

\bibitem{Guyon}
J.~Guyon.
\newblock Limit theorems for bifurcating {M}arkov chains. {A}pplication to the
  detection of cellular aging.
\newblock {\em Ann. Appl. Probab.}, 17(5-6):1538--1569, 2007.

\bibitem{gyo-pau2016}
B.~M. Gyori and D.~Paulin.
\newblock Hypothesis testing for markov chain monte carlo.
\newblock {\em Statistics and Computing}, 26(6):1281--1292, 2016.

\bibitem{hh:ml}
P.~Hall and C.~C. Heyde.
\newblock {\em Martingale limit theory and its application}.
\newblock Academic Press, Inc. [Harcourt Brace Jovanovich, Publishers], New
  York-London, 1980.
\newblock Probability and Mathematical Statistics.

\bibitem{Parzen1962}
E.~Parzen.
\newblock On estimation of a probability density function and mode.
\newblock {\em The Annals of Mathematical Statistics}, 33(3):1065--1076, 1962.

\bibitem{silverman1986density}
B.~W. Silverman.
\newblock {\em Density Estimation for Statistics and Data Analysis}, volume~26.
\newblock CRC Press, 1986.

\end{thebibliography}
\end{document}